\documentclass[12pt]{amsart}
\usepackage{amsfonts}
\usepackage{latexsym,amsmath,amsthm,amssymb,amsfonts,mathrsfs}
\usepackage{graphicx}
\numberwithin{equation}{section}

\newtheorem{thm}{Theorem}[section]

\newtheorem{lem}{Lemma}[section]
\newtheorem{prop}{Proposition}[section]
\newtheorem{cor}{Corollary}[section]
\theoremstyle{definition}
\newtheorem{dfn}{Definition}[section]

\theoremstyle{remark}
\newtheorem{rem}{Remark}

\newcommand{\il}{\langle}
\newcommand{\li}{\rangle}

\newcommand{\LL}{\bigg \langle}
\newcommand{\RR}{\bigg \rangle}

\begin{document}

\title[Lagrangian self-shrinkers and Piecewise mean curvature flow]{Geometry of Lagrangian self-shrinking tori and applications to the Piecewise Lagrangian Mean Curvature Flow}

\author{Jingyi Chen and John Man Shun Ma}

\address{Department of Mathematics,
The University of British Columbia, Vancouver, BC 
Canada V6T1Z2}
\email{jychen@math.ubc.edu.ca}
\email{johnma@math.ubc.edu.ca}

\begin{abstract}
We study geometric properties of the Lagrangian self-shrinking tori in $\mathbb R^4$.  When the area is bounded above uniformly, we prove 
that the entropy for the Lagrangian self-shrinking tori can only take finitely many values; this is done by deriving a {\L}ojasiewicz-Simon type gradient inequality for the branched conformal self-shrinking tori and then combining with the compactness theorem in \cite{CMa}.  When the area bound is small, we show that any Lagrangian self-shrinking torus in $\mathbb R^4$ with small area is embedded with uniform curvature estimates, and the space of such tori is compact.

Using the finiteness of entropy values, we construct a piecewise Lagrangian mean curvature flow for Lagrangian immersed tori in $\mathbb R^4$, along which the Lagrangian condition is preserved, area is decreasing, and the type I singularities that are compact with a fixed area upper bound can be perturbed away in finite steps. This is a Lagrangian version of the construction for embedded surfaces in $\mathbb R^3$  in \cite{CM} . 
\end{abstract}

\date{\today}

\maketitle

\section{Introduction}

One of the major challenging problems in the study of Lagrangian mean curvature flow is to formulate a weak version of the mean curvature flow that preserves the Lagrangian condition and goes beyond singular time, as  the well-known weak forms of mean curvature flow such as the Brakke flow or the level set approach do not work well in the Lagrangian setting. 

It is known that the rescaled mean curvature flow (MCF) at a finite time singularity converges to a self-shrinking solution, the so-called {\sl self-shrinker}; the convergence may be weak if the singularity is not of type I (cf. \cite{H2}, \cite{Ilmanen}, \cite{White}). As local models of singularities, self-shrinkers are important.  In \cite{CM}, Colding and Minicozzi introduced an entropy functional (see (\ref{definition of entropy})) of a hypersurface (cf. \cite{MM}) and showed that the sphere and the cylinders are the only entropy stable self-shrinking hypersurfaces. Using this and  a compactness theorem \cite{CM1} on the space of embedded self-shrinking surfaces in $\mathbb R^3$, they constructed in \cite{CM}  a piecewise MCF for embedded surfaces in $\mathbb R^3$ (under some assumptions), such that if a uniform diameter estimate holds 
 then the flow  shrinks to a round point. 

When the initial immersion is Lagrangian, a basic fact proved by Smoczyk \cite{Sm1} is that the MCF preserves the Lagrangian condition. 
In \cite{LL}, \cite{LZ} the authors studied the Lagrangian entropy stability of Lagrangian self-shrinking immersions and obtained entropy instability results. In particular, Li and Zhang showed in \cite{LZ} that if $F: M^n\to \mathbb R^{2n}$ is a closed orientable Lagrangian self-shrinker and the first Betti number of $M$ is greater than $1$ then $F$ is Lagrangian entropy unstable \footnote{More precisely, in \cite{LL}, \cite{LZ}, the authors study the Lagrangian $\mathcal F$-stability of a Lagrangian immersion. The equivalence of $\mathcal F$-stability and entropy stability is proved in \cite{CM} for the hypersurface case and can be generalized to higher codimensional cases. See section 2 and \cite{ALY}.}. Since there is no simply connected closed Lagrangian self-shrinker (see \cite{Sm0}, and \cite{CM} for branched immersions in dimension 2), all closed orientable Lagrangian self-shrinkers in $\mathbb R^4$ are Lagrangian entropy unstable. 

In \cite{CMa}, we used the harmonic map theory to show that if $\{ F_n :(\Sigma , h_n) \to \mathbb R^4\}$ is a sequence of compact Lagrangian self-shrinking immersions with uniform area bound and that the conformal structures $h_n$ do not degenerate, then a subsequence converges to a {\sl branched conformal} Lagrangian self-shrinker $F : (\Sigma, h)\to \mathbb R^4$. Furthermore, the assumption on the conformal structure can be dropped if $\Sigma$ is topologically a torus $\mathbb T$. This provides a compactification of the space of compact Lagrangian self-shrinking immersions in $\mathbb R^4$. A key ingredient in proving the compactness result is the rigidity established in \cite{CM}: there are no {\it branched} Lagrangian self-shrinking 2-spheres. This prevents development of the bubbles in Sacks-Uhlenbeck's compactness procedure \cite{SU1}, \cite{P}, hence yields strong convergence for the harmonic mappings.

The goal of this paper is twofold. First, we establish geometric properties of branched conformal Lagrangian self-shrinking tori. The discussion is divided into two cases:  for small area bound we prove embeddedness and curvature estimates, and for arbitrary area bound we show finite discreteness of the values of the entropy. To achieve the latter, we derive a {\L}ojasiewicz-Simon type
gradient inequality for the energy operator $\mathscr E$ naturally defined on the total space of $C^{2,\alpha}$-mappings and the moduli space of the conformal structures on the torus. This should have other applications. 
Second, the Lagrangian entropy instability, compactness for the space of Lagrangian self-shrinking tori and the finiteness of the value distribution of the entropy in Theorem \ref{finitely many entropy value on X_Lambda} below together lead us to define a piecewise Lagrangian MCF for a Lagrangian immersed torus in $\mathbb R^4$ which preserves the Lagrangian condition and the Maslov class, decreases area and avoids compact type I singularities with an arbitrarily given area upper bound in a finite number of steps. 

We now state our results on the Lagrangian self-shrinking tori. 
 
\begin{dfn} \label{definition of mathfrak X_Lambda}
Let $\Lambda$ be a positive number. Let $\mathfrak X_\Lambda$ be the space of branched conformally immersed Lagrangian self-shrinking tori with area less than or equal to $\Lambda$.
\end{dfn}

When the area upper bound $\Lambda$ is not small (as in Theorem \ref{Theorem for geometric property for small area}), it is not known whether any branched conformal Lagrangian self-shrinking torus with nonempty branch locus exists or not. The possible existence of branch points of elements in ${\mathfrak X}_\Lambda$  is a serious obstacle for applications to Lagrangian MCF as one would hope to perturb the branched Lagrangian surface to a nearby Lagrangian immersion, but such resolution of singularity in the Lagrangian setting, even in dimension two, is not available. 
Note that it is in general difficult to study nearby branched immersions by deforming them along the normal vector fields. In particular, it is hard to study stability problem of branched Lagrangian self-shrinking immersions as in \cite{CM}, \cite{LZ}, and Weinstein's Lagrangian neighbourhood theorem \cite{Weinstein} does not apply to the branched case. In view of all these and the special feature of the embedded graphic representation of a surface near a self-shrinker in the codimension one case, the idea of the piecewise MCF introduced in \cite{CM} is not directly applicable to the Lagrangian case in $\mathbb R^4$, even with the compactness theorems in \cite{CMa}. 

In order to construct a piecewise Lagrangian MCF for torus, we observe in this paper that one can bypass the issue of branchedness of a limiting surface in ${\mathfrak X}_{\Lambda}$ by controlling the entropy values $\lambda(F)$ attained by the self-shrinkers, where for a branched immersion $F:\mathbb T\to\mathbb R^4$ its entropy is defined by 
$$
\lambda(F) = \sup_{x_0\in\mathbb R^4,t_0>0}\frac{1}{4\pi t_0}\int_{\mathbb T} e^{-\frac{|F(x)-x_0|^2}{4t_0}}d\mu_F.
$$

The theorem below is a crucial ingredient in our construction of piecewise Lagrangian MCF for torus, but it is also interesting in its own right: it is equivalent to that in the induced metric from $G=e^{-\frac{|x|^2}{4}} \delta_{ij}$ on ${\mathbb R}^4$ the areas of branched Lagrangian self-shrinking tori in ${\mathfrak X}_{\Lambda}$ can only take a finite number (depending on $\Lambda$) of values for any given $\Lambda$.  

\begin{thm} \label{finitely many entropy value on X_Lambda}
Let $\lambda : \mathfrak X_{\Lambda} \to [0,\infty)$ be the entropy function which sends $F$ to its entropy $\lambda(F)$. Then the image of $\lambda$ is finite for any given $\Lambda$. 
\end{thm}

To prove Theorem \ref{finitely many entropy value on X_Lambda}, we derive a {\L}ojasiewicz-Simon gradient inequality for branched conformal self-shrinking 2-dimensional tori. The celebrated {\L}ojasiewicz-Simon gradient inequality is proved in \cite{Simon} with important applications to the harmonic map flow and the minimal cones. Since the pioneering work \cite{Simon}, the inequality and its variation has wide applications in geometric problems. For MCF, Schulze \cite{Schulze} used the inequality to prove a uniqueness result for compact embedded singularity of tangent flow. Colding and Minicozzi \cite{CM15} derived  {\L}ojasiewicz-Simon gradient inequalities in a noncompact setting and settled the uniqueness problem for all generic singularities of mean convex MCF at all singularities.

The classical {\L}ojasiewicz-Simon gradient inequality is established for real analytic functionals over a compact manifold whose Euler-Lagrange operator is elliptic and of order 2. In our case, we are concerned with the entropy functional $\lambda$, which is, at a self-shrinker, just the area of the shrinker in $(\mathbb R^4, G)$ up to a universal constant.  
However, in our situation, the self-shrinkers might be branched and the Euler-Lagrange operator of the area functional fails to be elliptic at the branch locus, so Simon's infinite dimensional version of the {\L}ojasiewicz inequality in \cite{Simon} is not directly applicable. To overcome the difficulty, we consider the real analytic energy functional $\mathscr E$  defined on the mapping space $C^{2,\alpha}(\mathbb T,\mathbb R^4)$ together with the Teichm\"uller space of ${\mathbb T}$, 
and continue to view self-shrinkers as branched minimal immersions in $(\mathbb R^4, G)$ \cite{Angenent}. 
The functional $\mathscr E$ has been extensively used in minimal surface theory, especially, in showing existence of minimal surfaces.
A critical point of $\mathscr E$ corresponds to a branched conformal self-shrinking torus. Since the space of conformal structures on a torus is two dimensional, the ellipticity of the $L^2$-gradient of $\mathscr E$ at a critical point of $\mathscr E$ for each {\it fixed} conformal structure enables us to show that the second order derivative $\mathcal L$ of $\mathscr E$ at the critical point is a Fredholm operator of index zero, which is sufficient to derive the desired gradient inequality. Theorem \ref{finitely many entropy value on X_Lambda} is then a direct consequence of the gradient inequality and the compactness results in \cite{CMa}.

If $\Lambda <32\pi$, the Willmore functional of a self-shrinker with area upper bound $\Lambda$ is less than $8\pi$; a classical theorem of Li and Yau \cite{LY} then asserts that all compact Lagrangian self-shrinking tori must be embedded and without branched point. Using recent results of Lamm-Sch\"atzle in \cite{LS} and the compactness results in \cite{CMa}, we show that the upper bound can be pushed beyond Li-Yau's estimate:   

\begin{thm} \label{Theorem for geometric property for small area}
There are positive numbers $\epsilon_0, \epsilon_1$ and $C_0$, where $\epsilon_1 \le \epsilon_0$, so that 
\begin{enumerate}
\item (Compactness) The space $\mathfrak X_{32\pi +\epsilon}$ is compact for all $0<\epsilon \le \epsilon_0$, 
\item (No Branch Points) All elements in $\mathfrak X_{32\pi+\epsilon_0}$ are immersed, and all elements in $ \mathfrak X_{32\pi+\epsilon_1}$ are embedded. 
\item (Curvature Estimates) If $F \in \mathfrak X_{32\pi+\epsilon_0}$, then the second fundamental form of $F$ is bounded by $C_0$. 
\end{enumerate}
\end{thm}

Lastly, we apply Theorem \ref{finitely many entropy value on X_Lambda} to construct a piecewise Lagrangian MCF for Lagrangian immersed torus $F :\mathbb T\to \mathbb R^4$ (see Definition \ref{definition of piecewise MCF}). In fact, we only need the result for compatified space of immersed Lagrangian self-shrinking tori. We show that all type I singularities with an arbitrarily given area upper bound can be perturbed in finitely many steps, where a smooth Lagrangian MCF for torus restarts at each step, such that the same kind of singularities will not appear in the last step. We remark that the perturbation can be made arbitrarily small while fixing the number of perturbation performed. Note that, in the special case of small area, Theorem \ref{Theorem for geometric property for small area} is sufficient 
since the existence of a nearby Lagrangian immersion of the torus around a limiting surface in $\mathfrak X_\Lambda$ (now immersed) follows from the Lagrangian neighbourhood theorem.  

\vspace{.2cm}

Our main result on Lagrangian mean curvature flow in a weak form is   

\begin{thm} \label{Generic LMCF}
Let $F : \mathbb T\to \mathbb R^4$ be an immersed Lagrangian torus and let $\Lambda, \delta >0$ be given constants. Then there exists a piecewise Lagrangian MCF $\{F^i_t : i = 0, 1, \cdots, k-1\}$ with initial condition $F$, where $k\le | \lambda (\mathfrak X_\Lambda)| <\infty$, such that the singularity at time $t_{k}$ is not a type I singularity modelled by a compact self-shrinker with area less than or equal to $\Lambda$. Moreover, the Maslov class of each immersion is invariant along the flow.
\end{thm}

Under an additional assumption, we prove a similar result in Theorem \ref{g>1} for the case of genus larger than one. 

The organization of this paper is as follows. Some background is provided in section 2. In section 3, we provide with proof necessary results in generalizing the compactness theorem in \cite{CMa} from immersions to branched immersions. The compactness result is stated in Theorem \ref{compactness for branch self-shrinkers}.  We then 
prove Theorem \ref{Theorem for geometric property for small area}. In section 4, we derive a {\L}ojasiewicz-Simon gradient inequality for branched conformal self-shrinking torus and prove Theorem \ref{finitely many entropy value on X_Lambda}. In the last section we prove Theorem \ref{Generic LMCF} and a result in the higher genus case. 

\vspace{.2cm}

{\bf Acknowledgement.} The first author is grateful for the partial support of an NSERC Discovery Grant (RGPIN 203199-1). Part of this work is supported by the National Science Foundation under Grant No. DMS-1440140 while the first author was in residence at the Mathematical Sciences Research Institute in Berkeley, California, during February 16 - March 31, 2016.


\section{Background}
\subsection{Mean curvature flow and the self-shrinkers}
A family of immersions $F_t : \Sigma\to \mathbb R^N$ from an $n$-dimensional manifold $\Sigma$ to the Euclidean space is said to satisfy the mean curvature flow (MCF) if
\begin{equation} \label{definition of MCF}
\frac{\partial F_t}{\partial t} =\vec H,
\end{equation}
where $\vec H$ is the mean curvature vector given by
\[
\vec H = \text{tr} A,\ \ \ A (X, Y) = (D_X Y)^\perp,
\]
where $A$ is the second fundamental form of the immersion and $\perp$ denotes the normal component of a vector. An immersion is called self-shrinking (or a self-shrinker) if it satisfies 
\begin{equation} \label{self-shrinking eqn}
\vec H = -\frac 12 F^\perp.
\end{equation}
If $F$ is self-shrinking, then up to a family of diffeomorphisms, the family of immersions 
$$\{\sqrt{-t} F: t\in [-1,0)\}$$
solves the MCF. The self-shrinkers model the singularity of MCF (cf. \cite{H2,Ilmanen, White}). 

An immersion $F: \Sigma \to \mathbb R^{2n}$ of an $n$-dimensional manifold $\Sigma$ is called Lagrangian if $F^*\omega = 0$, where 
\begin{equation} \label{definition of omega}
\omega = \sum_{i=1}^n dx^i \wedge dy^i 
\end{equation}
is the standard symplectic form on $\mathbb R^{2n}$. If $F_t :\Sigma\times[0,T)\to\mathbb R^{2n}$ is the MCF starting from a Lagrangian immersion $F_0$, then $F_t$ is Lagrangian for all $t \in [0,T)$ \cite{Sm1}. 

When $\Sigma$ is a surface, a branched conformal immersion $F : \Sigma\to \mathbb R^N$ is called a self-shrinker if (\ref{self-shrinking eqn}) is satisfied on $\Sigma\setminus B$, where $B$ is the set of branch points. A branched conformal immersion $F: \Sigma \to \mathbb R^4$ is called Lagrangian if $F^*\omega = 0$, where $\omega$ is as in (\ref{definition of omega}) (see \cite{CMa}, section 3). 

\subsection{Lagrangian $\mathcal F$-stability and Lagrangian entropy staibility}
The entropy $\lambda$ and $\mathcal F$-stability are introduced in \cite{CM} for an embedded self-shrinking hypersurfaces and are later carried over in \cite{ALY}, \cite{LL}, \cite{LZ} for all codimensions. The Lagrangian case is discussed in \cite{LL}, \cite{LZ} and the definition of Lagrangian $\mathcal F$-stability is introduced therein. We start with recalling the definitions of the $\mathcal F$ and $\lambda$ functionals and the related stability.   When we consider Lagrangian immersions, we will assume $N=2n$. 
\begin{dfn} Let $(x_0, t_0) \in \mathbb R^N \times \mathbb R_{>0}$. The $\mathcal F$-functional of an immersion $F : \Sigma^n \to \mathbb R^N$ is given by 
\begin{equation} \label{definition of mathcal F}
\mathcal F_{x_0, t_0} (F) = (4\pi t_0)^{-\frac n2} \int_\Sigma  e^{-\frac{|F(x)-x_0|^2}{4t_0}} d\mu_F.
\end{equation}
\end{dfn}

The $\mathcal F$-functional characterizes the self-shrinkers as follows:  $F:\Sigma \to \mathbb R^N$ is a self-shrinker if and only if 
\[ \frac{d}{ds} \mathcal F_{x_s, t_s} (F_s) \bigg|_{s=0} = 0\]
for all variations $(x_s, t_s, F_s)$ such that $(x_0, t_0, F_0) = (0,1, F)$. 

We recall that a normal vector field $X$ along a Lagrangian immersion is called a Lagrangian variation if 
\begin{equation}
 d( \iota_X \omega) = 0. 
\end{equation}


\begin{dfn}
A self-shrinker $F$ is called (Lagrangian) $\mathcal F$-stable if for all (Lagrangian) variations $F_s$, there is a variation $(x_s, t_s)$ so that 
\[
\frac{d^2}{ds^2} \mathcal F_{x_s, t_s} (F_s)\bigg|_{s=0} \ge 0.
\]
\end{dfn}

In \cite{LZ}, Li and Zhang calculated the second variation of the $\mathcal F$-functional of a Lagrangian immersion with respect to the Lagrangian variations. They proved:

\begin{thm} \label{Theorem of Li and Zhang}  
Let $\Sigma$ be a compact orientable $n$-dimensional manifold whose first Betti number is greater than $1$. If $F: \Sigma \to {\mathbb R}^{2n}$ is a Lagrangian self-shrinker, then $F$ is Lagrangian $\mathcal F$-unstable.
\end{thm} 

When $F : \Sigma \to \mathbb R^{2n}$ is a Lagrangian immersion, let $F_s:\Sigma\to\mathbb R^{2n}$ be a normal variation of $F$ such that each $F_s$ is a Lagrangian immersion. In this case, the normal variational vector field $X=\frac{d}{ds}|_{s=0}F_s$ can be identified with a closed 1-form on $\Sigma$ by $X\mapsto -\iota_X \omega$. The converse is also true as seen in the following elementary lemma. 

Recall that if $\alpha$ is a 1-form on a Riemannian manifold $(\Sigma, g)$ then $\alpha^\sharp$ is the vector field on $\Sigma$ uniquely determined by
\begin{equation} \label{definition of vector fields dual to 1-form alpha}
g(\alpha^\sharp, Y) = \alpha(Y), \ \ \ \forall \,Y \in T\Sigma.
\end{equation}

\begin{lem} \label{closed one form gives variations}
Let $F: \Sigma^n \to \mathbb R^{2n}$ be a Lagrangian immersion and let $\alpha$ be a closed 1-form on $\Sigma$. Then there is a family of Lagrangian immersions $F_s : \Sigma \to \mathbb R^{2n}$ so that $F_0 = F$ and
\begin{equation} \label{construct Lag variation using alpha}
\frac{d}{ds}\bigg|_{s=0} F_t = J \alpha^\sharp,
\end{equation}
where $J$ is the standard complex structure on $\mathbb R^{2n}$. 
\end{lem}

\begin{proof}
Let $\pi: N\Sigma\to \Sigma$ be the normal bundle of the immersion $F$. Then the mapping $$\tilde F(x, v) = F(x) + v$$ 
is a local diffeomorphism from a tubular neighbourhood $U$ of the zero section of $N\Sigma$ onto its image in $\mathbb R^{2n}$.

Since $\alpha$ is a closed 1-form on $\Sigma$, $\beta = (\pi|_U)^* \alpha$ is a closed 1-form on $U$, and $\beta$ sends the normal vectors $v$ to zero. The pullback 2-form $\omega_0 = \tilde F^* \omega$ on $U$ is closed as $\omega$ is closed and it is non-degenerate as $\tilde F$ is a locally diffeomorphic and $\omega$ is non-degenerate. Let $X$ be the vector field on $U$ dual to $\beta$ with respect to $\omega_0$, that is,
\begin{equation} \label{integrate alpha}
\beta (Y) = - \,\omega_0 ( X, Y)
\end{equation}
for all vector fields $Y$ on $U$. Let $\phi_s$ with $s\in (-\epsilon, \epsilon)$ be the one parameter group of diffeomorphisms on $U$ generated by $X$. Then $F_s := \tilde F\circ  \phi_s|_\Sigma : \Sigma \to \mathbb R^{2n}$ is a family of Lagrangian immersions in $\mathbb R^{2n}$ and $F_0 = \tilde F \circ \phi_0|_\Sigma = \tilde F|_\Sigma = F$.

It remains to verify (\ref{construct Lag variation using alpha}). By the definition of $\tilde{F}$, its differential $\tilde{F}_*$ maps the tangent vectors to the zero section $(\Sigma,0)$ at the point $(x,0)\in U$ to the tangent vectors to the image surface $F(\Sigma)$ at the point $F(x)\in{\mathbb R}^{2n}$ and it maps the normal vectors to the normal vectors by the identity map at the corresponding points. 
We need to check $X = J\alpha^\sharp$. Let $Y_1, Y_2$ be arbitrary  tangent vectors to the zero section $(\Sigma,0)$ at a point $(x,0)$. Since $JY_2$ is normal to $\Sigma$ as $\Sigma$ is Lagrangian and $\omega(X, Y) = \langle JX, Y\rangle$, we have
\[
\begin{split} \alpha(\pi_*Y_1)  &= \beta (Y_1 + JY_2) \\
&= -\,\omega_0(X, Y_1 + JY_2) \\
&= - \, \omega(\tilde{F}_*X, \tilde{F}_*Y_1+\tilde{F}_*JY_2)\\
&= -\, \omega(\tilde{F}_*X,\tilde{F}_*Y_1+JY_2)\\
&= -\,\il J\tilde{F}_*X, \tilde{F}_*Y_1+ JY_2\li \\
&= -\,\il J\tilde{F}_*X, \tilde{F}_*Y_1\li -\il \tilde{F}_*X, Y_2\li
\end{split}
\]
As $\tilde{F}$ is locally diffeomorphic, $X$ is normal to the zero section because $Y_2$ is arbitrary. Then it follows from the arbitrariness of $Y_1$ that $-JX = \alpha^\sharp$, by dropping the notion $\tilde{F}_*$. This is the same as $X = J\alpha^\sharp$.
\end{proof}

The entropy of a hypersurface is defined in \cite{CM, MM}. The definition for an immersion in any codimension is the same. 
\begin{dfn}
The entropy of an immersion $F: \Sigma \to \mathbb R^N$ is defined as
\begin{equation} \label{definition of entropy}
\lambda (F) = \sup_{x_0, t_0} \mathcal F_{x_0, t_0}(F).
\end{equation}
\end{dfn}
It is clear that $\lambda(F)$ is invariant under translations and scalings. Huisken's monotonicity formula \cite{H2} implies that $\lambda( F_t)$ is non-increasing if $\{F_t\}$ satisfies the MCF, and is constant if and only if $\{F_t\}$ is self-shrinking. Analogous to the entropy stability introduced in \cite{CM},  
we define Lagrangian entropy stability of a Lagrangian self-shrinker. 
\begin{dfn}
Let $F:\Sigma \to \mathbb R^{2n}$ be a self-shrinker. Then $F$ is called Lagrangian entropy stable if $\lambda(\tilde F) \ge \lambda (F)$ for all Lagrangian immersions $C^0$ close to $F$.
\end{dfn}
In \cite{CM}, it is proved that every $\mathcal F$-unstable embedded self-shrinking hypersurface which does not split off a line is entropy unstable. As observed in \cite{ALY}, the exact same proof works for any codimension. 
According to \cite{LZ}, the second variation formula for the $\mathcal F$-functional at a closed self-shrinker can be rewritten in terms of the closed 1-form dual to the Lagrangian variation field. Therefore, when $F:\Sigma \to \mathbb R^{2n}$ is a Lagrangian $\mathcal F$-unstable self-shrinker, there is a closed 1-form $\alpha$ on $\Sigma$ so that $\mathcal F''(\alpha) <0$ for all variations $(x_s, t_s)$ of $(0,1)$. To proceed from the Lagrangian $\mathcal F$-instability to the Lagrangian entropy instability, one needs to use the actual family $F_s$ of Lagrangian immersions coming from the Lagrangian variation. By Lemma \ref{closed one form gives variations}, there is a Lagrangian variation $\{F_s\}$ that corresponds to $\alpha$. By taking a family of diffeomorphism $\phi_s : \Sigma \to \Sigma$, we can further assume that $\{F_s\}$ is a family of normal variations. Thus the same proof of Theorem 0.15 in \cite{CM} can be carried over to show that $F: \Sigma \to \mathbb R^{2n}$ is also Lagrangian entropy unstable. We omit the proof here.

\begin{thm} \label{entropy unstable}
Let $\Sigma$ be compact and $F : \Sigma \to \mathbb R^{2n}$ be an immersed Lagrangian self-shrinker. If $F$ is Lagrangian $\mathcal F$-unstable, then it is also Lagrangian entropy unstable. In particular, there is a Lagrangian immersion $\widehat F :\Sigma \to \mathbb R^{2n}$ so that $\lambda (F) > \lambda(\widehat F)$. Moreover, $\widehat F$ can be chosen to be arbitrarily close to $F$, in the sense of smallness of $\|F - \widehat F\|_{C^k}$ for all $k$.
\end{thm}


\section{Lagrangian self-shrinking tori with small area}
In this section, we prove Theorem \ref{Theorem for geometric property for small area}. We will use a contradiction argument, and by doing so we need to extend the compactness theorems (Theorems 1.2 and 1.3 in \cite{CMa}) to {\it branched} conformally immersed Lagrangian self-shrinking surfaces. This extension will be done for any area upper bound (not necessarily small), and beside Theorem \ref{Theorem for geometric property for small area}, it will also be used in the next section in the proof of Theorem \ref{finitely many entropy value on X_Lambda}. Yet, for the construction of a piecewise Lagrangian MCF of Lagrangian torus, the compactness results in \cite{CMa} suffice.  

Then we will combine the compactness results with the factorization result of Lamm and Sch{\"a}tzle \cite{LS} concerning conformal immersion of torus into $\mathbb R^4$ with Willmore energy $8\pi$ to conclude the theorem. 

\begin{rem} \label{remark on different definition of self-shrinker}
In \cite{CMa}, a branched self-shrinker is defined as a branched immersion $F: \Sigma \to \mathbb R^N$ which satisfies 
\[
F^\perp = -\vec H,
\]
as opposed to (\ref{self-shrinking eqn}), which has an extra $1/2$ factor. Note that both definitions are common in the literatures and differ only by a scaling of the branched immersion. As a result, it should be pointed out that the constants that appear in this section are slightly different from those in \cite{CMa}. 
\end{rem}

As in \cite{CMa}, we view a branched self-shrinker $F:\Sigma \to \mathbb R^4$ in $\mathbb R^4$ as a harmonic map from $(\Sigma, h)$ to $(\mathbb R^4 , G)$. Here $h$ is the conformal structure on $\Sigma$ such that $F$ is conformal with respect to $h,$ and $G$, where $G$ is the metric on $\mathbb R^4$ given by 
\begin{equation} \label{definition of G}
G_{ij}= e^{-\frac{|x|^2}{4}} \delta_{ij},
\end{equation}
where $\delta_{ij}$ represents the standard Euclidean metric on $\mathbb R^4$. Then we use general harmonic map theories from \cite{P} and \cite{CLi}. In order to use these results, as in \cite{CMa}, we need to show that the self-shrinkers with a uniform area upper bound stay in a bounded domain in $(\mathbb R^4, G)$. In particular, we need the following lemma, which extends Lemma 4.1 in \cite{CMa} to allow branch points. The proof is almost the same as that of Lemma 4.1 in \cite{CMa}, except extra effort needs to be given at the branch points.

\begin{lem} \label{F_n in bounded ball}
Let $F$ be a compact branched conformal self-shrinker in $\mathbb R^4$. Then the image of $F$ lies in a ball of radius $R_0$ centered at the origin in $\mathbb R^4$, where $R_0$ depends only on $\mu(F)$.  
\end{lem}

\begin{proof}
Let $F : \Sigma \to \mathbb R^4$ be a branched conformally immersed self-shrinker. By (\ref{self-shrinking eqn}), the equation 
\begin{equation} \label{Lap |F|^2 1}
\Delta_g |F|^2 = -|F^\perp|^2 +4
\end{equation}
holds on $\Sigma\setminus B$, where $g = F^*\langle \cdot, \cdot\rangle$ and $B$ is the finite branch locus. 

First, we show that $F$ must intersect the closed ball centered at the origin of $\mathbb R^4$ with radius $2$. Since $F$ is a branched conformal immersion, there is a nonnegative smooth function $\varphi$ and a smooth metric $g_0$ on $\Sigma$ compatible with the conformal structure $h$ so that $g = \varphi g_0$. Therefore 
$$ 
\varphi \,\Delta_g = \Delta_{g_0}
$$ 
and by (\ref{Lap |F|^2 1}),
\begin{equation} \label{Lap |F|^2 g_0}
\Delta_{g_0} |F|^2 = \varphi \, (-|F^\perp|^2 +4).
\end{equation}

Unlike (\ref{Lap |F|^2 1}), (\ref{Lap |F|^2 g_0}) is satisfied everywhere on $\Sigma$, as both sides of the equation are continuous and $B$ is finite. Since $\Sigma$ is compact, the smooth function $|F|^2$ attains its minimum, say at $x_0\in \Sigma$. Since $F$ is a minimal immersion in $(\mathbb R^4, G)$, the tangential component $F^\top$ is well defined at a branch point and $F^\top (x_0) = 0$. If $F$ is immersed at $x_0$, by using (\ref{Lap |F|^2 1}) as in \cite{CMa}, the weak maximum principle shows that $|F(x_0)|^2 \le 4$ and we are done. Thus we only need to rule out the case that $F$ is branched at $x_0$, $|F(x_0)|^2 >4$ and there does not exist any immersed point $y \in \Sigma$ so that $|F(y)|^2 = |F(x_0)|^2$. Assume this case happens. Since the branch points are isolated, $|F|^2$ has a strict minimum at $x_0$. Noting that $F^\top (x_0) = 0$ and $|F|^2 = |F^\top|^2 + |F^\perp|^2$, we have $|F^\perp (x)|^2 >4$ in a neighbourhood of $x_0$. By (\ref{Lap |F|^2 g_0}) we have $\Delta_{g_0} |F|^2 \le 0$ in the neighbourhood. However, this contradicts the strong maximum principle, and we are done.

Next, we show that the extrinsic distance between any two points on the image of $F$ is bounded above by a constant that depends only on the area upper bound. 
Note that
\begin{equation} \label{Lap |F|^2 2}
\Delta_g |F|^2 d\mu_g = d*_g d|F|^2.
\end{equation}
and the Hodge star operator $*_g$ depends only on the conformal class of $g$, $\Delta_g |F|^2 d\mu_g$ is well-defined on $\Sigma$. Thus we integrate (\ref{Lap |F|^2 1}) and use (\ref{self-shrinking eqn}) to get
\begin{equation} \label{4W = area}
\mathcal W(F) := \frac 14 \int_\Sigma |\vec H|^2 d\mu = \frac 14 \mu(F).
\end{equation}
One also note that Simon's diameter estimate \cite{Simon1} holds for 2-varifolds with square integrable generalized mean curvature ((A.16) in \cite{KS}). Thus there is a constant $C$ such that
\[
\bigg(\frac{\mu (F)}{ \mathcal W(F)}\bigg)^{\frac{1}{2}} \leq  \text{diam}F(\Sigma) \leq C\big(  \mu (F) \mathcal W(F)\big)^{\frac{1}{2}} ,
\]
where
\[
\text{diam}F(\Sigma) := \sup_{x, y\in \Sigma} |F(x) - F(y)|.
\]
Together with (\ref{4W = area}), we see that 
$$
\text{diam} F(\Sigma)\leq \frac{1}{2}C\mu(F).
$$ 

It follows that the image of $F$ lies in $B(R_0)$ for some $R_0$ depending only on the area upper bound. 
\end{proof}

Let $\{F_n : \Sigma \to \mathbb R^4\}$ be a sequence of compact Lagrangian branched conformal self-shrinkers with uniform upper bound $\Lambda$. Lemma \ref{F_n in bounded ball} implies that the images of $F_n$ lie in a bounded region $B_R$ in $\mathbb R^4$. The Riemannian space $(B_R, G)$ can be isometrically embedded into a compact Riemannian manifold $(N,g)$, by Lemma 4.2 in \cite{CMa}. We can assume that $\{F_n\}$ is also a sequence of harmonic mappings from $(\Sigma, h_n)$ to the compact Riemannian manifold $(N, g)$. Using the same argument in \cite{CMa}, we extend Theorems 1.2 and 1.3 in \cite{CMa} to the branched immersions. 

\vspace{.2cm}

We now state the main compactness result.

\begin{thm} \label{compactness for branch self-shrinkers}
Let $F_n : \Sigma \to \mathbb R^4$ be a sequence of branched conformally immersed Lagrangian self-shrinkers with a uniform area upper bound $ \Lambda$.
\begin{enumerate}
\item If $\Sigma = \mathbb T$ is topologically a torus, then by passing to subsequence if necessary, there is a conformal structure $h$ on $\mathbb T$ so that $F_n : (\mathbb T, h_n) \to \mathbb R^4$ converge smoothly to a branched conformally immersed Lagrangian self-shrinker $F: (\mathbb T, h)\to \mathbb R^4$ and $h_n$ converge to $h$.
\item If $\Sigma$ is a closed surface of genus $g\ge 2$ and assume additionally that the conformal structures $h_n$ converge to a conformal structure $h$ on $\Sigma$. Then by passing to subsequence if necessary, $F_n : (\Sigma, h_n) \to \mathbb R^4$ converges smoothly to a branched conformally immersed Lagrangian self-shrinker $F: (\Sigma, h)\to \mathbb R^4$.
\end{enumerate}
\end{thm}

Note that in \cite{CMa} we first show (2) and then argue that the conformal structures must converge when $\Sigma = \mathbb T$ because degeneration of conformal structures in the genus 1 case would create Lagrangian self-shrinking 2-spheres (possibly branched) but this would violate our rigidity result (see \cite{CMa} for the details). 

In particular, (1) in Theorem \ref{compactness for branch self-shrinkers} implies (1) in Theorem \ref{Theorem for geometric property for small area}. The proof of the remaining parts of Theorem \ref{Theorem for geometric property for small area} will be divided into the following results. We recall that $\mathfrak X_\Lambda$ stands for the space of branched conformally immersed Lagrangian self-shrinking tori of area no larger than $\Lambda$.

\begin{prop} \label{no branched points in small area}
There is a positive number $\epsilon_0$ so that if $F \in \mathfrak X_{32\pi+\epsilon_0}$, then $F$ is immersed.
\end{prop}

\begin{proof}
Arguing by contradiction, we assume that there is a sequence $F_n : \mathbb T\to \mathbb R^4$ of branched conformal Lagnrangian self-shrinking tori so that
\begin{equation} \label{epsilon_0 gap equation 1}
 \liminf_{n\to \infty} \mu (F_n) \le  32\pi
\end{equation}
and each $F_n$ has a nonempty set of branch points. Using Theorem \ref{compactness for branch self-shrinkers}, by passing to a subsequence if necessary, the sequence $\{F_n\}$ converges smoothly to a branched conformal Lagrangian self-shrinking torus $F_\infty : \mathbb T\to \mathbb R^4$. Let $B_n$ be the set of branch points of $F_n$. Since $\mathbb T$ is compact, again by passing to a subsequence if necessary, there is a sequence $\{p_n\}$, where $p_n \in B_n$ for each $n\in \mathbb N$, so that $p_n \to p \in \mathbb T$. As $DF_n (p_n) = 0$ for all $n\in \mathbb N$ and the convergence $F_n \to F_\infty$ is smooth, $DF(p) = 0$ and so $p$ is a branch point of $F_\infty$, where $DF, DF_n$ are the differentials of $F,F_n$, respectively. By the theorem of Li and Yau (Theorem 6 in \cite{LY}, see also the appendix in \cite{KS} for the generalization to branched immersions), since $F_\infty$ is not embedded,
\begin{equation} \label{W(F) ge 8 pi}
\mathcal W (F_\infty) \ge 8\pi.
\end{equation}

On the other hand, from (\ref{epsilon_0 gap equation 1}) and Theorem 1 in \cite{CLi}, 
\[
\mu(F_\infty) \le \liminf \mu(F_n) \le 32\pi.
\]
Together with (\ref{W(F) ge 8 pi}) and (\ref{4W = area}) we have $\mathcal W(F_\infty) = 8\pi$. Since $\mathcal F_\infty$ has a branch point, Proposition 2.3 in \cite{LS} implies that $F_\infty$ factors through a branched conformal immersion $g: \mathbb T \to \mathbb S^2$. It follows that there is a branch conformal Lagrangian self-shrinking sphere $h: \mathbb S^2 \to \mathbb R^4$ so that $F_\infty = h\circ g$. However, by Theorem 1 in \cite{CMa}, such a $h$ does not exist. This contradicts the existence of the sequence $\{F_n\}$. The proposition is now proved. 
\end{proof}

Proposition \ref{no branched points in small area} and  Theorem \ref{compactness for branch self-shrinkers} lead to 

\begin{thm} \label{compactness theorem with small area}
Let $\epsilon_0$ be as in Proposition \ref{no branched points in small area}. Then the space of all Lagrangian {\bf immersed} self-shrinking tori with area less than or equal to $32 \pi + \epsilon_0$ is compact.
\end{thm}

Next we prove part (3) in Theorem \ref{Theorem for geometric property for small area}.
\begin{cor} \label{curvature estimates}
(Curvature Estimates) There is $C_0>0$ so that if $F: \mathbb T\to \mathbb R^4$ is a Lagrangian immersed self-shrinking torus with area less than or equals to $32\pi + \epsilon_0$, then the second fundamental form of $F$ is bounded by $C_0$.
\end{cor}

\begin{proof}
Assume this were not true. Then there is a sequence $F_n : \mathbb T \to \mathbb R^4$ of Lagrangian immersed self-shrinking tori with area less than $32\pi + \epsilon_0$ so that
\begin{equation}\label{Curvature estimates proof}
\max_{F_n(\mathbb T)}|A_n| \to \infty,
\end{equation}
where $A_n$ is the second fundamental form of the immersion $F_n$. Using Theorem \ref{compactness theorem with small area}, a subsequence of $\{F_n\}$ converges smoothly to an immersed self-shrinker $F_\infty$. In particular, we have
$$
(g_n)_{ij} = \frac{\partial F_n}{\partial x_i}\cdot \frac{\partial F_n}{\partial x_j} \longrightarrow  \frac{\partial F_\infty}{\partial x_i}\cdot \frac{\partial F_\infty}{\partial x_j} = (g_{\infty})_{ ij}, \,\,\, \mbox{as $n\to\infty$}.
$$
Since $g_\infty$ is positive definite as $F_\infty$ is immersed, there is a positive number $C$ so that $g_n \ge C \delta_{ij}$ for all $n$. So $g^{-1}_n$ are uniformly bounded. Hence 
$$
\max_{F_n(\mathbb T)}|A_n|^2 = \max_{F_n(\mathbb T)}g_n^{ij}g_n^{kl}\langle (A_n)_{ik},(A_n)_{jl}\rangle 
$$ 
are uniformly bounded and (\ref{Curvature estimates proof}) is impossible. 
\end{proof}

To finish the proof of Theorem \ref{Theorem for geometric property for small area}, it remains to prove the second part in (2).

\begin{prop}
There is a positive constant $\epsilon_1 \le \epsilon_0$ so that if $F \in \mathfrak X_{32\pi +\epsilon_1}$, then $F$ is embedded. 
\end{prop} 

\begin{proof}
As in the proof of Corollary \ref{curvature estimates}, assume the contrary that there is a sequence $\{F_n\}$ of immersed, non-embedded Lagrangian self-shrinking tori with $\mu(F_n)\le 32\pi +\epsilon_0$ and $\mu(F_n) \to 32\pi$. By Theorem \ref{compactness theorem with small area}, after passing to a subsequence if necessary, $\{F_n\}$ converges smoothly to an immersed Lagrangian self-shrinking torus $F_\infty : \mathbb T \to \mathbb R^4$ with area $\mu(F_\infty) = 32\pi$. By (\ref{4W = area}), the Willmore energy of $F_\infty$ is $8\pi$. Since each  $F_n$ is non-embedded, there are distinct points $p_n, q_n \in \mathbb T$ so that
\begin{equation} \label{F_n (p_n) = F_n(q_n)}
F_n (p_n) = F_n(q_n).
\end{equation}
As $\mathbb T$ is compact, we may assume $p_n \to p$ and $q_n \to q$. Taking $n\to \infty$ in (\ref{F_n (p_n) = F_n(q_n)}), we have $F_\infty(p) = F_\infty (q)$. First of all, we must have $p=q$: Indeed, if $p\neq q$, then $F_\infty$ is not embedded and that contradicts Theorem 2.2 in \cite{LS}, which states that any immersion $F :\mathbb T \to \mathbb R^4$ with $\mathcal W(F) = 8\pi$ has to be embedded.

Let $d_n$ be the distance function on $\mathbb T$ induced by the pullback metric $F_n^*\langle \cdot, \cdot \rangle$. As $p = q$ and $\{F_n\}$ converges smoothly to $F$, we have $\ell_n : = d_n (p_n, q_n) \to 0$ as $n\to \infty$. Let $\eta_n : [0,\ell_n] \to \mathbb T^2$ be a shortest geodesics in $(\mathbb T, F_n^*\langle\cdot, \cdot\rangle)$ joining $p_n$ to $q_n$. Since $F_n(\eta_n(0)) = F_n(\eta_n(\ell_n))$, $F_n \circ \eta_n : [0, \ell_n] \to \mathbb R^4$ is a closed curve in $\mathbb R^4$ with length $\ell_n$. Let $\gamma_n : [0, \ell_n] \to \mathbb R^4$ be the translation $\gamma_n (t) = F_n\circ \eta_n (t) - F_n(p_n)$. Then each $\gamma_n$ is parameterized by arc length and $\gamma_n(0) = \gamma_n(\ell_n) = 0 \in \mathbb R^4$. Using the following simple estimates
\[
\begin{split}
\ell_n &= \int_0^{\ell_n} \langle \gamma'_n(t) , \gamma'_n(t) \rangle dt\\
&= -\int_0^{\ell_n} \langle \gamma_n (t), \gamma_n''(t)\rangle dt + \langle \gamma_n(\ell_n), \gamma'_n(d_n)\rangle - \langle \gamma_n (0), \gamma'_n(0)\rangle \\
&= -\int_0^{\ell_n} \langle \gamma_n (t), \gamma''_n (t)\rangle dt \\
&\le \int_0^{\ell_n} | \gamma_n (t)| \cdot | \gamma''_n(t)| dt \\
&\le \ell_n \int_0^{\ell_n} |\gamma''_n (t)| dt,
\end{split}
\]
we obtain
$$\int_0^{\ell_n} |\gamma''_n (t) | \,dt \ge 1.$$
Since $\ell_n \to 0$, the above inequality implies that there is $s_n \in [0, \ell_n]$ so that $| \gamma''_n (s_n)| \to \infty$ as $n\to \infty$. Since $\eta_n$ is a geodesic on $(\mathbb T, F_n^*\langle\cdot, \cdot\rangle)$,
$$
\gamma''_n = (F_n\circ \eta_n)'' = \nabla^n_{\eta_n'} \eta_n' + A_n(\eta_n', \eta'_n) = A_n (\eta_n', \eta_n'),
$$
where $\nabla^n$ is the Levi-Civita connection on $(\mathbb T, F^*_n\langle \cdot, \cdot\rangle)$ and $A_n$ is the second fundamental form of $F_n(\mathbb T)$ in ${\mathbb R}^4$. Thus 
$$
| \gamma_n''(t)| \le |A_n (\eta_n(t))|
$$ 
and this implies
$$
\max_{F_n(\mathbb T)} | A_n |^2 \to \infty
$$
as $n\to \infty$. However, this is impossible by Corollary \ref{curvature estimates}.
\end{proof}


\section{A {\L}ojasiewicz-Simon type gradient inequality for branched self-shrinking tori}
In the last section we show that with a small area bound, all Lagrangian self-shrinking tori are immersed at least. This makes it much easier to study the space $\mathfrak X_{32\pi+\epsilon_0}$, as all nearby Lagrangian self-shrinking tori can be deformed to each other by using the normal vectors fields. However, it is difficult in general to relate two nearby branched conformal immersions, even if they are $C^k$-close when treated as mappings to the Euclidean space. In particular, it seems difficult to extend the perturbation procedure as in \cite{CM, LZ}, where the stability condition is described by using the normal vector fields, to {\it branched} conformal self-shrinkers, .

In this section, we show that the entropy $\lambda$ is  locally a constant function in the space of branched conformal compact self-shrinking tori $F :\mathbb T \to \mathbb R^4$. To do this we derive a {\L}ojasiewicz-Simon type gradient inequality for branched conformal self-shrinking tori $F: \mathbb T \to \mathbb R^4$. In the genus one case, the explicit expression of the conformal structures in the Teichm\"uller space makes the computation and the real analyticity of the functional ${\mathscr E}$ transparent.  
Once this is done, together with the compactness of $\mathfrak X_\Lambda$, we conclude the proof of Theorem \ref{finitely many entropy value on X_Lambda}. 

\subsection{A Fredholm operator of index zero}  Let $(\Sigma, g)$ be a compact Riemannian surface and $(M, h)$ a Riemannian manifold. Given a $C^1$ mapping $F: \Sigma \to M$, the energy of $F$ is given by 
\[ \mathscr E_{g,h} (F)= \frac 12 \int_\Sigma e_{g,h} (F) d\mu_g,\]
where $e_{g,h}(F)$ is the norm of the differential $DF _x : T_x\Sigma \to T_{F(x)}M$. Locally it is given by
\[
e_{g,h} (F) =  g^{ij} h_{\alpha\beta}\frac{\partial F^\alpha}{\partial x^i} \frac{\partial F^\beta}{\partial x^j} .
\]

For a fixed $h$, define $\mathscr E: C^1(\Sigma,M)\times \{g: \mbox{$g$ is a Riemannian metric on $\Sigma$}\} \to\mathbb R$ by 
$$
\mathscr E(F,g)=\mathscr E_{g,h}(F).
$$  

\begin{lem} \label{conformal implies critical wrt g}
 If $F : (\Sigma, g)\to (M, h)$ is conformal, then $g$ is a critical point of $\mathscr E$ with respect to all its smooth variations $g_s$, where $g_0=g$. That is, 
\[ \frac{d}{ds}{\mathscr E}_{g_s, h} (F) \bigg|_{s=0} = 0. \]
\end{lem}

\begin{proof}
Let $g_s$ be a family of smooth metrics on $\Sigma$ so that $g_0 = g$ and $\dot g = \frac{d}{ds} g_s\big|_{s=0}$. Then
\[
\begin{split}
\frac{d}{ds}\left(g^{ij} \sqrt{\det g}\right)\bigg|_{s=0}  &= - g^{ik}g^{jl} \dot g_{kl} \sqrt{\det g} + \frac 12 g^{ij}\sqrt{\det g} g^{kl}\dot g_{kl}\\
&= \left( \frac 12 g^{kl}g^{ij} - g^{ik} g^{jl}\right) \dot g_{kl} \sqrt{\det g}
\end{split}
\]
Thus
\begin{equation}  \label{calculation of d E wrt g}
\frac{d}{ds}e_{g_s, h} (F) d\mu_g \bigg|_{s=0} = \left( \frac 12 g^{kl}g^{ij} - g^{ik} g^{jl}\right) \dot g_{kl}h_{\alpha\beta}\frac{\partial F^\alpha}{\partial x^i} \frac{\partial F^\beta}{\partial x^j} d\mu_g
\end{equation}
Since $F$ is conformal,
\begin{equation} \label{F is conformal}
h_{\alpha\beta}\frac{\partial F^\alpha}{\partial x^i} \frac{\partial F^\beta}{\partial x^j} = \varphi g_{ij}
\end{equation}
for some function $\varphi$ on $\Sigma$. Put (\ref{F is conformal}) into (\ref{calculation of d E wrt g}) and use $g^{ij} g_{ij} = 2$ since $\Sigma$ is two dimensional, we see that $\frac{d}{ds}\mathscr E_{g_s, h} (F)\big|_{s=0}  =0$, as claimed.
\end{proof}
On the other hand, recall that a branched minimal immersion is (weakly) conformal and harmonic, and we have the following (\cite{SU1}, Theorem 1.8)
\begin{prop} \label{critical wrt g implies conformal}
If $u$ is critical map of $\mathscr E$ with respect to the variations of $u$ and the conformal structures on $\Sigma$, then $u$ is a branched minimal immersion.
\end{prop}

Let $ U$ be an open subset in the upper half space $\mathbb H = \{ \tau\in \mathbb C \, | \,\text{Im} \,\tau >0\}$. It is well-known that the upper half space represents the Teichm{\" u}ller space of the standard torus $\mathbb T = \mathbb R^2/\{1,i\}$ and we treat $U$ as a local parameterization of the conformal structures on $\mathbb T$ near a given one.

Let $0<\alpha <1$ be fixed. Define
\[
\begin{split}
\mathscr U &= C^{2, \alpha} (\mathbb T, \mathbb R^4) \times  U, \\
\mathscr C^{k,\alpha}&= C^{k,\alpha}(\mathbb T, \mathbb R^4) \oplus \mathbb R^2 \\
\mathscr W^{k,p} &= W^{k,p} (\mathbb T, \mathbb R^4) \oplus \mathbb R^2 \\
\mathscr L^2 &=\mathscr W^{0,2}.
\end{split}
\]
Note that $\mathscr C^{k,\alpha}, \mathscr W^{k,\alpha}$ are Banach spaces \footnote{All Banach spaces considered in this paper are real Banach spaces.} with the norms
\begin{eqnarray*}
&&\| (\phi, \nu) \|_{k,\alpha} = \|\phi \|_{C^{k, \alpha}} + |v|,  \\
&&\|(\phi, \nu)\|_{\mathscr W^{k,p}} = \|\phi\|_{W^{k,p}} + |v|
\end{eqnarray*}
respectively. When $(M,h)=({\mathbb R}^4,G)$, where $G$ is as in (\ref{definition of G}), the functional $\mathscr E : \mathscr U \to \mathbb R$ takes the form 
\begin{equation} \label{definition of mathscr E}
\mathscr E(u, \tau) =\frac 12 \int_{\mathbb T} e^{-\frac{|u|^2}{4}} |Du|_{\tau}^2d\mu _\tau,\ \ \ (u,\tau) \in \mathscr U. 
\end{equation}
Here $g_\tau$ is the metric on $\mathbb T$ given by
\begin{equation} \label{definition of g_tau}
g_\tau  =\left( \begin{matrix} 1 & \tau_1 \\ 0  &\tau_2 \end{matrix}\right)^T\left( \begin{matrix} 1 & \tau_1 \\ 0  &\tau_2 \end{matrix}\right)
\end{equation}
and
\[
d\mu_\tau = d\mu_{g_\tau} = \sqrt{\det g_\tau} \, dxdy, \ \ \ |Du|^2_\tau = g^{ij}_\tau D_i u \cdot D_j  u.
\]
The metric $g_\tau$ is in the conformal class represented by $\tau$, as it can be seen easily that $g_\tau$ is the pullback metric via the linear mapping from $\mathbb T={\mathbb R}^2/\{1,i\}$ to ${\mathbb R}^2/\{1,\tau\}$. 
Note that for each fixed $\tau$, $\mathscr E(\cdot, \tau)$ is the Dirichlet energy functional of the mappings $u:(\mathbb T, g_\tau) \to (\mathbb R^4, G)$.

It is well-known \cite{Angenent} that minimal surface in $(\mathbb R^4, G)$ corresponds to self-shrinking surfaces in $\mathbb R^4$. Thus Lemma \ref{conformal implies critical wrt g} and Proposition \ref{critical wrt g implies conformal} imply the following

\begin{prop} \label{harmonic and conformal iff critical pt of E}
$(u, \tau )$ is a critical point of $\mathscr E$ if and only if $u :(\mathbb T, g_\tau)\to \mathbb R^4$ is a branched conformal self-shrinking torus.
\end{prop}

Next we consider the $L^2$-gradient $\mathscr M:\mathscr U \to \mathscr C^{0,\alpha}$ of $\mathscr E$. That is, we find for each $(u, \tau) \in \mathscr U$ an element $\mathscr M(u, \tau) \in \mathscr C^{0,\alpha}$ so that for all $(\phi, \nu)\in \mathscr C^{2, \alpha}$,

\begin{equation} \label{definition of mathscr M}
\frac{d}{ds}\bigg|_{s=0} \mathscr E(u + s \phi, \tau + s\nu) = \langle \mathscr M(u, \tau), (\phi, \nu)\rangle_{u, \tau}.
\end{equation}
Here we define
\[
\langle \phi_1, \phi_2\rangle_{u,\tau}  =\int_\mathbb{T} \phi_1 \cdot \phi_2 \;\mathrm e^{-\frac{|u|^2}{4}} d\mu_\tau
\]
and
\begin{equation} \label{definition of <>_{u,tau}}
 \langle (\phi_1, \nu_1), (\phi_2, \nu_2)\rangle_{u, \tau} = \langle \phi_1, \phi_2\rangle_{u,\tau}+ \nu_1 \cdot \nu_2.
\end{equation}
Integrating by parts, we see that

\begin{equation} \label{expression of mathscr M}
\mathscr M(u, \tau) = \left(-g^{ij}_\tau e^{\frac{|u|^2}{4}} D_j (e^{-\frac{|u|^2}{4}} D_i u ) - \frac 14|Du|^2_\tau u , \nabla \mathscr E^u_\tau\right)
\end{equation}
where $\mathscr E^u: U \to \mathbb R$ is given by $\mathscr E^u (\tau) = \mathscr E(u, \tau)$ and $\nabla \mathscr E^u_\tau$ is the gradient of $\mathscr E^u$ at $\tau$.  

 Let $(u, \tau)$ be a critical point of $\mathscr E$, that is, $\mathscr M (u, \tau)  =0$. Let 
$$
\mathcal L = \mathcal L_{(u, \tau)}:\mathscr C^{2,\alpha} \to \mathscr C^{0,\alpha}
$$ 
be the Fr{\'e}chet derivative of $\mathscr M$ at $(u, \tau)$. We will show that
\begin{equation} \label{expression of mathcal L}
\mathcal L(\phi, \nu) = \left(L \phi + \nabla_\nu B , (\nabla^2 \mathscr E^u_\tau) \nu + \langle \nabla B_\tau , \phi\rangle_{u, \tau} \right)
\end{equation}
where
\begin{equation}\label{L}
 L \phi= -g^{ij}_\tau e^{\frac{|u|^2}{4}} D_j (e^{-\frac{|u|^2}{4}} D_i\phi) - \frac 14 |Du|^2_\tau \phi + \frac 12g^{ij}_\tau D_j(u\cdot \phi)D_i u - \frac 12 g^{ij}_\tau (D_j u\cdot D_i \phi)u
\end{equation}
and $\nabla^2\mathscr E^u_\tau$ is the Hessian of $ \mathscr E^u$ at $\tau$; furthermore, $B: U\to C^{0,\alpha}(\mathbb T, \mathbb R^4)$ is given by
\begin{equation} \label{expression of B}
B(\sigma) =-g^{ij}_{\sigma} \left(  e^{\frac{|u|^2}{4}} D_j (e^{-\frac{|u|^2}{4}} D_iu) +\frac 14 (D_i u \cdot D_j u) u\right)
\end{equation}
and $\nabla B_\tau$ denotes the Fr{\'e}chet derivative of $B$ at $\tau$ and $\nabla_\nu B_\tau$ stands for the Fr{\'e}chet derivative of $B$ at $\tau$ in the direction $\nu$:
$$
\nabla_\nu B_\tau = \frac{d}{ds}\bigg|_{s=0} B(\tau+s\nu).
$$ 

To derive (\ref{expression of mathcal L}), note that the two terms in the first component of (\ref{expression of mathcal L}) arise from  direct differentiation of the first component of (\ref{expression of mathscr M}) with respect to $\phi$ and $\nu$. To derive the second component, note that $(\nabla^2\mathscr E^u_\tau) \nu$ is just the directional derivative of $\nabla \mathscr E^u_\tau$ with respect to $\nu$. Thus we need to show that $\nabla_\phi \nabla \mathscr E^u_\tau = \langle \nabla B_\tau, \phi\rangle_{u,\tau}$, where 
\[
\nabla_\phi  \nabla \mathscr E^u_\tau = \frac{d}{ds}\bigg|_{s=0} \nabla \mathscr E^{u+s\phi}_\tau. 
\]
 Note 
\[
\nabla\mathscr E^u_\tau = \frac 12\int_\mathbb{T} (\nabla g^{ij}_\tau) e^{-\frac{|u|^2}{4}} (D_i u \cdot D_j u) d\mu_\tau + \frac 14 \text{tr} (g^{-1}_\tau \nabla g_\tau)\mathscr E(u,\tau),
\]
where the second term on the right comes from differentiating the volume form $d\mu_\tau$. Since $(u, \tau)$ is a critical point of $\mathscr E$, this term vanishes when we differentiate with respect to $\phi$. Using this observation and integration by parts, 
\[
\begin{split}
\nabla_\phi \nabla \mathscr E^u_\tau &= \frac 12\int_\mathbb{T} (\nabla g^{ij}_\tau) \nabla_\phi \left(e^{-\frac{|u|^2}{4}} (D_i u \cdot D_j u)\right)d\mu_\tau \\
&= \frac 12 \int_\mathbb{T} (\nabla g^{ij}_\tau)  \left(-\frac 12 (u\cdot \phi) e^{-\frac{|u|^2}{4}} (D_i u \cdot D_j u)+2 e^{-\frac{|u|^2}{4}} (D_i u\cdot D_j\phi) \right) d\mu_\tau \\
&= - \int_\mathbb{T} (\nabla g^{ij}_\tau) \left(\frac 14 (D_i u \cdot D_j u ) u+ e^{\frac{|u|^2}{4}} D_j (e^{-\frac{|u|^2}{4}} D_i) \cdot \phi\right) e^{-\frac{|u|^2}{4}} d\mu_\tau \\
&= \int_\mathbb{T} \nabla B_\tau \cdot \phi \, e^{-\frac{|u|^2}{4}} d\mu_\tau \\
&= \langle \nabla B_\tau, \phi \rangle _{u,\tau}.
\end{split}
\]
Thus (\ref{expression of mathcal L}) is shown. 

\begin{lem}
Let $(u,\tau)$ be a critical point of ${\mathscr E}$. 
For all $(\phi,\nu), (\psi, \eta) \in \mathscr C^{2,\alpha}$, we have 
\begin{equation} \label{mathcal L is self adjoint}
\langle \mathcal L(\phi, \nu), (\psi, \eta)\rangle_{u,\tau} = \langle (\phi, \nu), \mathcal L (\psi, \eta)\rangle_{u,\tau}.
\end{equation}
\end{lem}

\begin{proof}
Let $\phi, \psi \in C^{2,\alpha}(\mathbb T, \mathbb R^4)$, then from (\ref{expression of mathcal L}) and (\ref{expression of B}),
\begin{equation} \label{L self-adjoint first line}
\begin{split}
\langle  L \phi, \psi\rangle_{u,\tau} &= \langle g^{ij}_\tau D_i \phi , D_j \psi\rangle_{u,\tau} -\frac 14 \langle |Du|^2 _\tau \phi, \psi\rangle_{u,\tau} \\
&+ \frac 12\int g^{ij}_\tau D_j(u \cdot \phi) D_i u \cdot \psi \,e^{-\frac{|u|^2}{4}} d\mu_\tau - \frac 12 \int g^{ij}_\tau (D_j u \cdot D_i\phi) (u\cdot \psi) \,e^{-\frac{|u|^2}{4}} d\mu_\tau.
\end{split}
\end{equation}
Integrating by parts for the third term on the right hand side in (\ref{L self-adjoint first line}) gives
\begin{equation} \label{L self-adjoint second line}
\begin{split}
\frac 12 \int g^{ij}_\tau D_j&(u \cdot \phi) D_i u \cdot \psi \,e^{-\frac{|u|^2}{4}} d\mu_\tau \\
&= -\frac 12\int g^{ij}_\tau (u\cdot \phi)(D_i u\cdot D_j \psi) \,e^{-\frac{|u|^2}{4}} d\mu_\tau - \frac 12\int (u\cdot \phi)g^{ij}_\tau  D_j (e^{-\frac{|u|^2}{4}} D_i u) \cdot \psi \,d\mu_\tau.
\end{split}
\end{equation}
Since $\mathscr M(u, \tau) = 0$, we have by (\ref{expression of mathscr M})
\[
g^{ij}_\tau D_j (e^{-\frac{|u|^2}{4}} D_i u) = -\frac 14 e^{-\frac{|u|^2}{4}} |Du|^2_\tau u.
\]
Putting this into (\ref{L self-adjoint second line}), we have
\[
\begin{split}
\langle L \phi, \psi\rangle_{u,\tau} &= \langle g^{ij}_\tau D_i \phi , D_j \psi\rangle_{u,\tau} -\frac 14 \langle |Du|^2 _\tau \phi, \psi\rangle_{u,\tau} \\
&\ \ -\frac 12\int g^{ij}_\tau (u\cdot \phi)(D_i u\cdot D_j \psi) e^{-\frac{|u|^2}{4}} d\mu_{g_\tau} -\frac 12\int g^{ij}_\tau (D_j u \cdot D_i\phi) (u\cdot \psi) e^{-\frac{|u|^2}{4}} d\mu_\tau\\
&\ \  + \frac 18 \int (u\cdot \phi)(u\cdot \psi) |Du|^2_\tau e^{-\frac{|u|^2}{4}} d\mu_\tau.
\end{split}
\]
Note that the right hand side is symmetric in $\phi$ and $\psi$. Thus
\begin{equation} \label{L is self adjoint}
\langle  L \phi, \psi\rangle_{u,\tau} = \langle \phi, L \psi\rangle _{u,\tau},\ \ \ \forall \phi, \psi\in C^{2, \alpha}(\mathbb T, \mathbb R^4).
\end{equation}
Using this, we have
\[
\begin{split}
\langle \mathcal L(\phi, \nu), (\psi, \eta)\rangle_{u,\tau} &= \langle L \phi + \nabla_\nu B_\tau , \psi\rangle_{u,\tau} + (\nabla^2 \mathscr E^u_\tau \nu + \langle \nabla B_\tau, \psi\rangle_{u,\tau}) \cdot  \eta \\
&= \langle  L \phi, \psi\rangle_{u,\tau} + \langle \nabla _\nu B_\tau , \psi\rangle_{u,\tau} + \langle\nabla_\eta B_\tau , \phi\rangle_{u,\tau} + (\nabla^2 \mathscr E^u_\tau \nu)\cdot \eta
\end{split}
\]
Again, the right hand side is symmetric in $(\phi, \nu)$ and $(\psi, \eta)$. We can now conclude the proof of the lemma.
\end{proof}

\begin{rem} Note that the apparent self-adjointness expression for $\mathcal L$ in (\ref{mathcal L is self adjoint}) only holds in ${\mathscr C}^{2,\alpha}$, and $\mathcal L$ is an operator from ${\mathscr C}^{2,\alpha}$ to ${\mathscr C}^{0,\alpha}$. Nevertheless, (\ref{mathcal L is self adjoint}) is useful in proving the following theorem. 
\end{rem}

\begin{thm}\label{L is Fredholm of index zero}
$\mathcal L$ is a Fredholm operator of index zero at a critical point $(u,\tau)$ of ${\mathscr E}$. 
\end{thm}

\begin{proof} The proof will be divided into several steps. 

{\bf Step 1.} We show that $\dim\ker\mathcal L$ is finite. 

Consider the first component of $\mathcal L$,
\begin{equation} \label{first component of mathcal L = 0}
L \phi + \nabla_\nu B_\tau  = 0.
\end{equation}
This equation is bilinear in $\phi,\nu$. Let $S$ be the subspace of $\mathbb R^2$ so that $\nu \in S$ if and only if (\ref{first component of mathcal L = 0}) has a solution. If $S = \{(0,0)\}$, then $\text{dim}\ker \mathcal L = \text{dim}\ker L<\infty$ since $L$ is elliptic. If not, let $\{\nu_i\}$ be a basis of $S$. Pick $\phi_i \in C^{2,\alpha}(\mathbb T, \mathbb R^4)$ so that $\phi_i$ satisfies (\ref{first component of mathcal L = 0}) with $\nu = \nu_i$. Let $(\phi, \nu)\in \ker\mathcal L$. Then $\nu \in S$. Write $\nu = \sum_i s^i \nu_i$ for some $s^i \in \mathbb R$. Then $\phi - s^i \phi_i \in \ker L$ and thus
\[ (\phi, \nu) = (\phi_0 , 0 ) + \sum_i s^i (\phi_i , \nu_i)\]
for some $\phi_0 \in \ker L$. Again, due to the ellipticity of $L$, $\dim\ker L$ is finite, hence $\ker \mathcal L$ is finite dimensional.

{\bf Step 2.} 
$\mathcal L$ has finite dimensional cokernel. Moreover, $\dim\ker{\mathcal L}=\dim  \text{coker}\mathcal L$.

 We will show that the mapping
\begin{equation} \label{mapping between ker mathcal L and coker mathcal L}
\ker \mathcal L \hookrightarrow \mathscr C^{2,\alpha}\hookrightarrow \mathscr C^{0, \alpha} \overset{\pi}{\to} \text{coker}\mathcal L
\end{equation}
is bijective, where $\pi$ is the projection to the quotient $\text{coker}\mathcal L = \mathscr C^{0,\alpha}/\text{Im}\mathcal L$. 

Firstly, if $(\psi_1, \eta_1), (\psi_2, \eta_2) \in \ker\mathcal L$ represent the same element in $\text{coker} L$, then there is $(\phi, \nu) \in \mathscr C^{2, \alpha}$ so that 
$$(\psi, \eta):= (\psi_1 - \psi_2, \eta_1 - \eta_2) = \mathcal L(\phi, \nu).$$ Using (\ref{mathcal L is self adjoint}),
\[
\langle(\psi, \eta), (\psi, \eta)\rangle_{u, \tau} = \langle \mathcal L(\phi, \nu), (\psi, \eta)\rangle_{u, \tau} = \langle (\phi, \nu), \mathcal L(\psi, \eta)\rangle_{u, \tau} = 0.
\]
Thus $(\psi, \eta) = 0$ and so the mapping $\ker\mathcal L \to \text{coker} \mathcal L$ defined in (\ref{mapping between ker mathcal L and coker mathcal L}) is injective.

Secondly,  we show that the mapping $\ker\mathcal L \to \text{coker} \,\mathcal L$ is surjective. 
Let $\overline{\text{Im}\,\mathcal L}$ be the $L^2$ closure of the image of $\mathcal L$ in $\mathscr L^2$ with respect to the inner product defined in (\ref{definition of <>_{u,tau}}). Let $(\psi, \eta) \in \mathscr C^{0, \alpha}$ represents an element in $\text{coker} \,\mathcal L$. We decompose $(\psi , \eta)$ into the component in $\overline{\text{Im}\,\mathcal L}$ and $\overline{\text{Im}\,\mathcal L}^\perp$. That is,
\begin{equation} \label{decomposition of (psi, eta)}
(\psi, \eta) = (\psi^\top, \eta^\top) + (\psi ^\perp, \eta^\perp)
\end{equation}
for some $\psi^\top, \psi^\perp \in L^2(\mathbb T, \mathbb R^4)$. Note that 
$$
\langle (\psi^\perp, \eta^\perp), \mathcal L(\phi,\nu)\rangle_{u,\tau}= 0
$$ 
for all $(\phi, \nu) \in \mathscr C^{2, \alpha}$. Letting $\nu =0$ and using (\ref{expression of mathcal L}), we have
\[
\langle L \phi , \psi^\perp \rangle_{u,\tau}  +\langle\nabla_{\eta^\perp} B_\tau, \phi\rangle_{u,\tau}  = 0,\ \ \ \forall \phi \in C^{2,\alpha}(\mathbb T, \mathbb R^4). 
\]
Note that the above equation is of the form
\[
\int (-g^{ij}_\tau D_{ij}\phi + \mathcal A^i D_i \phi + \mathcal B\phi)\cdot \psi^\perp dxdy= \int \mathcal F\cdot \phi \,dxdy,
\]
where $\mathcal A^i=(\mathcal A^i_{\beta \gamma})$ and $ \mathcal B = (\mathcal B_{\beta\gamma})$ are $(4\times 4)$-matrix-valued smooth functions and $\mathcal F = (\mathcal F_\beta)$ is a $\mathbb R^4$-valued smooth function. If we choose $\phi = (\rho, 0, 0,0)$, where $\rho \in C^\infty(\mathbb T, \mathbb R)$, we have
\begin{equation} \label{to use elliptic reg. of folland}
\int (- g^{ij}_\tau D_{ij} \rho+\mathcal A^i_{11} D_i \rho + \mathcal B_{11} \rho)\psi^\perp_1 dxdy = D(\rho),
\end{equation}
where
\begin{equation} \label{ distribution D}
D(\rho) = -\int  \sum_{k\neq 1} \mathcal A^i_{k1} \psi^\perp_k D_i \rho \,dxdy -\int \sum_{k\neq 1} \mathcal B_{k1} \psi^\perp_k \rho \,dxdy + \int \mathcal F_1 \rho \,dxdy.
\end{equation}
Since $\psi^\perp_k$ are in $L^2$ (noting that the $L^2$ spaces with respect the area elements $e^{-\frac{|u|^2}{4}}d\mu_\tau$ and $dxdy$ coincide over $\mathbb T$), as a distribution, $D$ is in $H^{-1}_{loc}$. Thus the Elliptic Regularity Theorem (Theorem 6.33 in \cite{Folland}) asserts $\psi_1^\perp \in H^1_{loc}$. Similarly, we have $\psi^{\perp}_k \in H^1_{loc}$ for $k=2,3,4$. Putting this information into (\ref{ distribution D}), we see that $D \in H^0_{loc}$, and in turn, this implies $\phi^\perp_1 \in H^2_{loc}$ by the Elliptic Regularity Theorem again. By a standard bootstrapping argument and the Sobolev embedding theorem, we see that $\psi^\perp \in C^{2,\alpha}$ (in fact, smooth). Using (\ref{mathcal L is self adjoint}) and the definition of $(\psi^\perp, \eta^\perp)$, we have 
$$(\mathcal L(\psi^\perp, \eta^\perp), (\phi,\nu))_{u,\tau} = 0$$
 for all $(\phi,\nu)\in{\mathscr C}^{2,\alpha}$, thus 
$$
\mathcal L(\psi^\perp, \eta^\perp) = 0.
$$

The smoothness of $(\psi^\perp,\eta^\perp)$ asserts $(\psi^\top,\eta^\top)\in{\mathscr C}^{0,\alpha}$. 
If we can show that 
\begin{equation} \label{Im mathcal L is closed in L^2 norm}
(\psi^\top, \eta^\top) \in \text{Im}\,\mathcal L,
\end{equation}
then $\pi(\psi, \eta) = \pi(\psi^\perp, \eta^\perp)$ by (\ref{decomposition of (psi, eta)}) and it follows that the mapping $\ker \mathcal L \to \text{coker} \mathcal L$ defined in (\ref{mapping between ker mathcal L and coker mathcal L}) is  surjective and we are done. To show (\ref{Im mathcal L is closed in L^2 norm}), recall that  $(\psi^\top, \eta^\top)\in \overline{\text{Im}\,\mathcal L}$. Thus there is a sequence $(\phi_n, \nu_n) \in \mathscr C^{2, \alpha}$ so that $\mathcal L(\phi_n, \nu_n) \to (\psi^\top, \eta^\top)$ in $\mathscr L^2$. Using the $\mathscr L^2$ inner product, we decompose $(\phi_n, \nu_n)$ into  
\begin{equation} \label{decomposition of (phi_n, nu_n)}
(\phi_n, \nu_n) = ( \phi_n^{\text K}, \nu_n^{\text K}) + (\phi_n^{\text P}, \nu_n^{\text P}),
\end{equation}
where $( \phi_n^{\text K}, \nu_n^{\text K}) \in \ker \mathcal L$ and $ (\phi_n^{\text P}, \nu_n^{\text P}) \in \ker \mathcal L^\perp$. Then by setting
\begin{equation} \label{definition of (psi_n, eta_n)}
(\psi_n, \eta_n) = \mathcal L (\phi^{\text P}_n, \nu^{\text P}_n)
\end{equation}
and using $\mathcal L(\phi_n^{\text K}, \nu_n^{\text K}) = 0$, we have
\[\begin{split}
(\psi_n, \eta_n) &= \mathcal L(\phi_n^{\text P}, \nu_n^{\text P})\\
&= \mathcal L(\phi_n, \nu_n) \\
&\overset{\mathscr L^2}{\to} (\psi^\top, \eta^\top) .
\end{split}
\]
The convergence above in particular implies that $\|\psi_n\|_{L^2} \le C$ for some constant $C$. From the first component of (\ref{expression of mathcal L}), which is 
\[ L\phi_n^{\text P} = \psi_n - \nabla _{\nu^{\text P}_n} B_\tau,\]
the standard elliptic estimates (Theorem 9.11 in \cite{GT}) implies that there are constants $C', C'', C''' >0$ so that 
\begin{equation} \label{L^2 estimate for phi_n^perp}
\begin{split}
\|\phi_n^{\text P}\|_{W^{2,2}} &\le C' \left(\|\phi_n^{\text P}\|_{L^2} + \| \psi_n - \nabla _{\nu_n^{\text P}} B_\tau\|_{L^2}\right) \\
&\le C'\left(\|\phi_n^{\text P}\|_{L^2} + \| \psi_n \|_{L^2} + C''|\nu_n^{\text P}|\right) \\
&\le C''' \left(\| (\phi_n^{\text P}, \nu_n^{\text P})\|_{\mathscr L^2} + 1\right). 
\end{split}
\end{equation}

Next, we show that the sequence $\{\|(\phi_n^{\text P}, \nu_n^{\text P})\|_{\mathscr L^2}\}$ is bounded. Assume not, then by taking a subsequence if necessary, we have $\|(\phi_n^{\text P}, \nu_n^{\text P})\|_{\mathscr L^{2}} \to \infty$. Let 
\begin{equation} \label{definition of (tilde phi_n, tilde nu_n)}
(\tilde \phi_n, \tilde \nu_n ) =\frac{  (\phi_n^{\text P}, \nu_n^{\text P})}{\|(\phi_n^{\text P}, \nu_n^{\text P})\|_{\mathscr L^2}}.
\end{equation}
Then, as $(\psi_n,\eta_n)$ converges to $(\psi^\top,\eta^\top)$ in ${\mathscr L}^2$,   
\begin{equation} \label{L(tilde phi_n, tilde nu_n) to 0}
\mathcal L(\tilde \phi_n, \tilde \nu_n) = \frac{(\psi_n, \eta_n)}{\|(\phi_n^{\text P}, \nu_n^{\text P})\|_{\mathscr L^2}} \overset{\mathscr L^2}{\longrightarrow } 0. 
\end{equation}
Since $\|( \tilde \phi_n, \tilde \nu_n)\|_{\mathscr L^2} = 1$, we may assume $\tilde \nu_n \to \tilde \nu$ for some $\tilde \nu \in \mathbb R^2$. From (\ref{L^2 estimate for phi_n^perp}) and (\ref{definition of (tilde phi_n, tilde nu_n)}), the sequence $\{ \|\tilde \phi_n\|_{W^{2,2}}\}$ is bounded. Hence, again by taking subsequence if necessary, there is $\tilde \phi \in W^{2,2}(\mathbb T, \mathbb R^4)$ so that $\tilde \phi_n \to \tilde \phi$ in $W^{1,2}(\mathbb T, \mathbb R^4)$. Using (\ref{L(tilde phi_n, tilde nu_n) to 0}), we have 
\[
\begin{cases}
L\tilde \phi + \nabla_{\tilde\nu} B_\tau = 0 & \text{weakly in } W^{1,2}(\mathbb T, \mathbb R^4), \\
\nabla^2\mathscr E^u_\tau \tilde \nu + \langle \nabla B_\tau, \tilde \phi\rangle_{u,\tau} = 0
\end{cases}
\]
Since $\tilde \phi \in W^{2,2}(\mathbb T,\mathbb R^4)$, the first equation is actually satisfied strongly in $W^{2,2}(\mathbb T, \mathbb R^4)$. Since $\nabla_{\tilde \nu} B_\tau$ is smooth, by the elliptic regularity, $\tilde \phi$ is smooth. Thus $(\tilde \phi, \tilde \nu) \in \mathscr C^{2,\alpha}$ and $\mathcal L(\tilde \phi, \tilde\nu) = 0$, in other words,
$(\tilde \phi, \tilde \nu) \in \ker\mathcal L .$
On the other hand, since $(\tilde \phi_n , \tilde \nu_n) \to (\tilde \phi, \tilde \nu)$ in $\mathscr L^2$ and $(\tilde \phi_n , \tilde \nu_n) \in \ker\mathcal L^\perp$, we also have $(\tilde \phi, \tilde \nu) \in \ker \mathcal L^\perp.$
Thus $(\tilde \phi, \tilde \nu) = (0,0)$. But this is impossible as $\|(\tilde \phi, \tilde \nu)\|_{\mathscr L^2}=1$ since $(\tilde \phi_n , \tilde \nu_n) \to (\tilde \phi, \tilde \nu)$ in $\mathscr L^2$ and $\|(\tilde \phi_n, \tilde\nu_n)\|_{\mathscr L^2} = 1$. The contradiction leads to the conclusion that the sequence $\{ \|(\phi_n^{\text P}, \nu_n^{\text P})\|_{\mathscr L^2}\}$ is bounded. 

From (\ref{L^2 estimate for phi_n^perp}), the sequence $\{ \|(\phi_n^{\text P}, \nu_n^{\text P})\|_{\mathscr W^{2,2}}\}$ is also bounded. By taking a subsequence if necessary, there is $(\phi, \nu) \in \mathscr W^{2,2}$ so that $(\phi_n^{\text P}, \nu_n^{\text P}) \to (\phi, \nu)$ in $\mathscr W^{1,2}$ and 
$$\mathcal L(\phi, \nu) = (\psi^\top, \nu^\top).$$ 
The first component of this is given by 
$$L \phi + \nabla_\nu B_\tau = \psi^\top.$$
Since $\phi \in W^{2,2}(\mathbb T, \mathbb R^4)$ and $\psi^\top \in C^{0, \alpha}(\mathbb T, \mathbb R^4)$, the standard elliptic regularity (Theorem 9.19 in \cite{GT}) implies that $\phi \in C^{2, \alpha}(\mathbb T, \mathbb R^4)$. Thus $(\phi, \nu) \in \mathscr C^{2, \alpha}$. This shows $(\psi^\top, \eta^\top) \in \text{Im}\mathcal L$. Therefore, the mapping $\ker\mathcal L\to\text{coker}\mathcal L$ is surjective.

{\bf Step 3.}  From the previous two steps, the bounded operator $\mathcal L$ has finite dimensional kernel and cokernel so it is a Fredholm operator of 
$$
\text{index}\,\mathcal L= \dim\ker\mathcal L - \dim\text{coker}\mathcal L=0.
$$
 This completes the proof of the theorem. 
\end{proof}

\subsection{A {\L}ojasiewicz-Simon type inequality} 
Next we prove a {\L}ojasiewicz-Simon gradient inequality for compact branched self-shrinkers $F: \mathbb T \to \mathbb R^4$. As in \cite{Simon}, we use the Liapunov-Schmidt reduction argument and the classical {\L}ojasiewicz inequality in \cite{Loja}. See \cite{FM} for a {\L}ojasiewicz-Simon inequality in the abstract setting and the related work in the reference therein. 

Let $$\Pi: \mathscr{L}^{2}\to\ker \mathcal L$$ be the $L^2$-projection with respect to the $L^2$ inner product:
\begin{equation} \label{standard mathscr L^2 inner product}
\langle (\psi_1, \nu_1), (\psi_2, \nu_2)\rangle_{\mathscr L^2} = \int_\mathbb{T} \psi_1\cdot\psi_2 \,dxdy + \nu_1 \cdot \nu_2
\end{equation}
for all $(\psi_1,\nu_1),(\psi_2,\nu_2)\in{\mathscr L}^{2}$. Recall that $\ker \mathcal L$ is a finite dimensional subspace and $\ker\mathcal L \subset \mathscr C^{\infty}$. For all $k= 0,1,2,\cdots$, we let 
$$\Pi_{k} : \mathscr C^{k,\alpha} \to \mathscr C^{0,\alpha}$$
 be the restriction of $\Pi$ to $\mathscr C^{k,\alpha}$ composed with the inclusion $\ker \mathcal L \hookrightarrow \mathscr C^{0,\alpha}$.  

\begin{lem} \label{Pi is bounded on C^k,alpha}
$\Pi_{k} :\mathscr C^{k, \alpha} \to \mathscr C^{0,\alpha}$ is a bounded linear operator for all nonnegative integers $k$. In particular, there is a positive constant $C_\alpha$ so that 
\begin{equation} \label{Pi is bounded equation}
\| \Pi_{k} (\psi, \nu)\|_{0, \alpha} \le  C_\alpha \| (\psi, \nu)\|_{k, \alpha}
\end{equation}
for all $(\psi, \nu) \in \mathscr C^{k, \alpha}$. 
\end{lem}

\begin{proof}
Let $(\chi_1, \nu_1), \cdots, (\chi_n ,\nu_n) \in \ker\mathcal L$ be an orthonormal basis of the finite dimensional space $\ker\mathcal L$ with respect to the inner product in (\ref{standard mathscr L^2 inner product}). Then for any $(\psi, \nu) \in \mathscr L^2$, we have
\[ \Pi (\psi, \nu) =\sum_{i=1}^n \langle (\chi_i, \nu_i), (\psi, \nu) \rangle_{ \mathscr L^2} (\chi_i, \nu_i).\]
Then we have
\[\begin{split}
\| \Pi_{k} (\psi, \nu)\|_{0, \alpha} &\le \sum^n_{i=1} \left| \langle (\chi_i, \nu_i), (\psi, \nu)\rangle_{ \mathscr L^2}\right| \| (\chi_i, \nu_i)\|_{0, \alpha} \\
&\le \left(\sum^n_{i=1} \| (\chi_i, \nu_i)\|_{0, \alpha} \right) \| (\psi, \nu)\|_{\mathscr L^2}.
\end{split}\]
Note that we used the Cauchy-Schwarz inequality and that $\| (\chi_i, \nu_i)\|_{\mathscr L^2} = 1$. Since 
$$
\int_{\mathbb T} dxdy  =1,
$$
we have   
\begin{equation} \label{L^2 to C^k,alpha is bounded with norm one}
\|(\psi, \nu)\|_{\mathscr L^2} \le \max_{\mathbb T^2}|\psi| +|\nu| \le  \| (\psi, \nu)\|_{k,\alpha}
\end{equation}
 for all nonnegative $k$. Now (\ref{Pi is bounded equation}) follows with $C_\alpha= \sum^n_{i=1} \| (\chi_i, \nu_i)\|_{0, \alpha}$.
\end{proof}

To simplify notations, in the sequel we use $x, y$ and $a, b$ to denote elements in $\mathscr C^{2, \alpha}$ and $\mathscr C^{0,\alpha}$ respectively. Let $x_{c} = (u, \tau)$ be a critical point of $\mathscr E$ as before, that is $\mathscr M(x_c) = 0$.

Consider the mapping $\mathscr N : \mathscr U \to \mathscr C^{0, \alpha}$ given by
\begin{equation} \label{definition of Phi}
\mathscr N (x) =\mathscr M (x)+\Pi_{2}(x-x_c). 
\end{equation}
Since $\Pi_2$ is linear, the differential $D\mathscr N$ at $x_c$ is given by
\begin{equation} \label{D Phi_{x_c} = mathcal L +Pi}
D\mathscr N_{x_c} = \mathcal L + \Pi_{2}.
\end{equation}
 
\begin{lem}
$D\mathscr N_{x_c}$ is bijective and its inverse is bounded. 
\end{lem}

\begin{proof}
First we show that $D\mathscr N_{x_c}$ is injective. Let $D\mathscr N_{x_c} (x) = 0$. Then by (\ref{D Phi_{x_c} = mathcal L +Pi}) we have 
$$\mathcal L(x) = -\Pi_{2} x.$$ Using (\ref{mathcal L is self adjoint}), for all $y\in \ker\mathcal L$ we have 
\[
\langle \Pi_{2} x, y\rangle_{u, \tau} = - \langle\mathcal L x , y\rangle_{u, \tau} =-\langle x, \mathcal L y\rangle_{u, \tau} = 0. 
\]
This means that $\Pi _{2} x\in\ker\mathcal L$ is orthogonal to $\ker\mathcal L$. Therefore, $\Pi_{2} x = 0$. Thus $\mathcal L x = 0$ and so $x \in \ker\mathcal L$. Hence $x = \Pi_{2} x = 0$ and $D\mathscr N_{x_c}$ is injective.

By Theorem \ref{L is Fredholm of index zero}, $\mathcal L$ is a Fredholm operator of index zero. Since $\Pi _{2}$ is bounded with a finite dimensional range, $\Pi_{2}$ is a compact operator and $D\mathscr N_{x_c} :\mathscr C^{2,\alpha} \to \mathscr C^{0, \alpha}$ is Fredholm with index zero (Theorem 5.10 in \cite{Schechter}). Together with the fact that $D\mathscr N_{x_c}$ is injective, $D\mathscr N_{x_c}$ is also surjective. Finally, the bounded inverse theorem (Theorem 3.8 in \cite{Schechter}) asserts that $D\mathscr N_{x_c}$ has a bounded inverse.
\end{proof}

By the inverse function theorem for Banach spaces (Theorem 15.2  in \cite{Deimling}), since $\mathscr N$ is $C^1$ ($\mathscr N$ is even analytic: see the appendix), there are open neighbourhoods $\mathscr U_1$ of $x_c$ in $\mathscr U$ and $\mathscr V_1$ of $0$ in $\mathscr C^{0, \alpha}$ so that $\mathscr N : \mathscr U_1 \to \mathscr V_1$ is invertible with a $C^1$ inverse $\Psi$. By shrinking $\mathscr U_1, \mathscr V_1$ if necessary, we assume that $\mathscr V_1$ is convex, $\mathscr U_1$ is contained in a convex set $\mathscr U_2 \subset\mathscr U$ and (since $\mathscr M$ and $\Psi$ are $C^1$) there exist two positive constants $M_1, M_2$ so that 
\begin{equation} \label{bounds on DPhi, DPsi}
\begin{split}
\| D\Psi(a)\|_{op} &\le M_1,\ \ \ \forall a\in \mathscr V_1, \\
\| D\mathscr M (x)\|_{op} &\le M_2,\ \ \ \forall x\in \mathscr U_2,
\end{split}
\end{equation}
where $\|\cdot\|_{op}$ denotes the operator norm for the corresponding operator. Using the Fundamental Theorem of Calculus, the above imply
\begin{equation} \label{Lipschitzness of Psi}
\| \Psi (a) -\Psi (b)\|_{2, \alpha} \le M_1 \| a-b\|_{0, \alpha}
\end{equation}
for all $a, b\in \mathscr V_1$ and 
\begin{equation} \label{Lipschitzness of M}
\| \mathscr M (x) -\mathscr M(y)\|_{0, \alpha} \le M_2 \| x-y\|_{2, \alpha}
\end{equation}
for all $x, y\in \mathscr U_1$. 

A main technical result in this section is the following {\L}ojasiewicz-Simon type gradient inequality:
\begin{thm} \label{prop: S L inequality for self-shrinking torus}
There is an open neighbourhood $\mathscr W_0 \subset \mathscr U$ of $x_c$, a positive constant $C_2$  and a constant $\theta \in (0,1/2)$ depending on $\mathscr E$ and $x_c$ so that
\begin{equation} \label{S L inequality}
|\mathscr E (x) - \mathscr E(x_c)|^{1-\theta} \le C_2 \|\mathscr M (x)\|_{0,\alpha},\ \ \ \forall x\in \mathscr W_0.
\end{equation}
\end{thm}

\begin{proof}
 Since $\Pi_{0}$ is bounded, there is an open neighbourhood $\mathscr V_0$ of $0$ so that $\mathscr V_0, \Pi_{0} \mathscr V_0 \subseteq\mathscr V_1$. For all $a\in \mathscr V_0$, $\Pi_0 a\in \mathscr V_1$. Since $\mathscr U_2$ is convex, the line segment joining $\Psi(a)$ and $\Psi(\Pi_0 a)$ is in $\mathscr U_2$. The Fundamental Theorem of Calculus and (\ref{definition of mathscr M}) yield
\[
\begin{split}
\mathscr E(\Psi(a)) - \mathscr E (\Psi(\Pi_{0} a)) &= -\int_0^1 \frac{d}{dt} (\mathscr E (\Psi(a) + t(\Psi(\Pi_{0} a) - \Psi(a))) dt \\
&= - \int_0^1 \langle \mathscr M (\Psi(a) + t(\Psi(\Pi_{0} a) - \Psi(a))), \Psi(\Pi_{0} a) - \Psi(a)\rangle_{u_t,\tau_t} dt,
\end{split}
\]
where we write 
$$ (u_t, \tau_t) = \Psi(a) + t(\Psi(\Pi_0 a) - \Psi(a)).$$
Using the Cauchy-Schwarz inequality, (\ref{L^2 to C^k,alpha is bounded with norm one}), (\ref{Lipschitzness of M}) and $|t|\le 1$, 
\begin{equation} \label{SL inequality prop first eqn}
\begin{split}
|\mathscr E&(\Psi(a)) -  \mathscr E (\Psi(\Pi_0 a))|\\
 & \le \|\mathscr M (\Psi(a) + t(\Psi(\Pi_{0}a) - \Psi(a)))\|_{ \mathscr L^2}\| \Psi(\Pi_{0} a) - \Psi(a)\|_{ \mathscr L^2} \\
&\le \|\mathscr M (\Psi(a) + t(\Psi(\Pi_{0} a) - \Psi(a)))\|_{0,\alpha} \| \Psi(\Pi_{0} a) - \Psi(a)\|_{2,\alpha}\\
&\le \big( \|\mathscr M (\Psi(a) \|_{0,\alpha} + M_2 \, t\, \|\Psi(\Pi_{0} a) - \Psi(a)\|_{2, \alpha}\big)\| \Psi(\Pi_{0} a) - \Psi(a)\|_{2,\alpha}\\
&\le \big( \|\mathscr M (\Psi(a) \|_{0,\alpha} + M_2\|\Psi(\Pi_{0} a) - \Psi(a)\|_{2, \alpha}\big)\| \Psi(\Pi_0 a) - \Psi(a)\|_{2,\alpha}
\end{split}
\end{equation}
On the order hand, since $a, \Pi_0a \in \mathscr V_1$, by (\ref{Lipschitzness of Psi}) we have 
\begin{equation} \label{Lemma 2.10 second eqn}
\| \Psi(\Pi_{0} a)-\Psi(a)\|_{2,\alpha} \le M_1 \| \Pi_{0} a- a\|_{0,\alpha}.
\end{equation}
Using the definition of $\mathscr N, \Psi$ and $\Pi_{0} \Pi_{2} = \Pi_{2}$, 
\begin{equation}\label{trivial}
a= \mathscr N(\Psi(a)) = \mathscr{M}(\Psi(a))+\Pi_2(\Psi(a)-x_c)
\end{equation}
\begin{equation} \label{Lemma 2.10 third eqn}
\begin{split}
\Pi_{0} a-a &= \Pi_{0} a - \mathscr M(\Psi (a)) - \Pi_{2} (\Psi(a) - x_c) \\
&= \Pi_{0} \big(a - \Pi_{2} (\Psi(a) - x_c)\big) - \mathscr M \big(\Psi(a)\big).
\end{split}
\end{equation}
Since $\Pi_0$ is bounded by Lemma \ref{Pi is bounded on C^k,alpha},
\[
\begin{split}
\|\Pi_{0} \big(a - \Pi_{2} (\Psi(a) - x_c)\big)\|_{0,\alpha} &\le C_\alpha \| a - \Pi_{2} (\Psi(a) - x_c)\|_{0,\alpha} \\
&= C_\alpha \|\mathscr M(\Psi(a))\|_{0,\alpha},
\end{split}
\]
where in the last line we use \eqref{trivial} again. Combining this with (\ref{Lemma 2.10 second eqn}) and (\ref{Lemma 2.10 third eqn}), we are led to 
\begin{equation} \label{S L inequality lemma equation}
\| \Psi (\Pi_{0} a) - \Psi (a)\|_{2, \alpha} \le C_1\| \mathscr M(\Psi(a))\|_{0,\alpha} 
\end{equation} 
for all $a\in \mathscr V_0$ with $C_1 = M_1(C_\alpha +1)$. Putting this into (\ref{SL inequality prop first eqn}), we have 
\begin{equation} \label{E(Psi (a)) - E(Psi (Pi a)) < C mathscr M^2}
|\mathscr E(\Psi(a)) -  \mathscr E (\Psi(\Pi_0 a))| \le C_3\| \mathscr M(\Psi(a))\|^2_{0,\alpha}
\end{equation}
for all $a\in \mathscr V_0$ and for some $C_3 >0$. 

Let $f : \mathscr V_1 \cap \ker\mathcal L \to \mathbb R$ be defined by 
\begin{equation} \label{definition of Gamma}
f (a) = \mathscr E (\Psi(a)).
\end{equation}
It is easy to show that $\mathscr E$, $\mathscr M$ are analytic (a proof is given in the appendix for completeness). Since $\Pi_{2}$ is linear, 
$$
\mathscr N = \mathscr M + \Pi_{2}-\Pi_2(x_c)
$$
 is analytic as well. Hence $\Psi$ is analytic by the analytic version of inverse function theorem (Theorem 15.3 in \cite{Deimling}). Consequently, as a composition of analytic functions, $f$ is also analytic, and it is defined on an open set in $\ker\mathcal L$, which is finite dimensional. The classical {\L}ojasiewicz inequality \cite{Loja} then implies that there is an open neighbourhood $\mathscr V_2 \subset \mathscr V_0$, constants $c >0$ and $\theta \in (0,1/2)$ so that 
\begin{equation} \label{classical Lojasiewicz inequality}
|f(\xi) - f (0)|^{1-\theta} \le c | f'(\xi)|,\ \ \ \forall \xi\in \mathscr V_2\cap \ker \mathcal L. 
\end{equation}
Using (\ref{definition of Gamma}) and (\ref{definition of mathscr M}), for all $b\in \mathscr V_1\cap K$ we have 
$$f'(b)(\cdot) = \langle \mathscr M(\Psi(b)) , D\Psi_b(\cdot)\rangle _{u, \tau}.$$
Using (\ref{L^2 to C^k,alpha is bounded with norm one}), (\ref{Lipschitzness of M}) and (\ref{S L inequality lemma equation}),
\begin{equation} \label{f' (Pi a) < C mathscr M (Psi a)}
\begin{split}
|f'(\Pi_{0} a)|& \le M_1 \|\mathscr M(\Psi (\Pi_{0} a)) \|_{\mathscr L^2} \\
&\le M_1 \| \mathscr M ( \Psi(\Pi_{0} a))\|_{0,\alpha}\\
&\le M_1 \big( \| \mathscr M ( \Psi(\Pi_{0} a)) - \mathscr M(\Psi (a))\|_{0,\alpha} + \|\mathscr M (\Psi(a)) \| _{0,\alpha} \big) \\
&\le M_1 \big( M_2 \| \Psi(\Pi_0 a) - \Psi (a)\|_{2,\alpha} +\|\mathscr M (\Psi(a)) \| _{0,\alpha}\big) \\
&\le C_4 \|\mathscr M (\Psi(a))\|_{0,\alpha}
\end{split}
\end{equation}
for some $C_4>0$. Now let $\mathscr W_0 = \Psi(\mathscr V_2)$. Thus for every $x\in \mathscr W_0$, there exists an $a\in \mathscr V_2$ such that $x = \Psi(a)$. By (\ref{f' (Pi a) < C mathscr M (Psi a)}), the classical {\L}ojasiewicz inequality (\ref{classical Lojasiewicz inequality}) and (\ref{E(Psi (a)) - E(Psi (Pi a)) < C mathscr M^2}),
\begin{equation}
\begin{split}
C_4 c\|\mathscr M(x)\|_{0,\alpha} &\ge c|f'(\Pi_0 a) |\\
&\ge |f(\Pi_0a) - f(0) |^{1-\theta} \\
& = |\mathscr E(\Psi( \Pi_0a)) -\mathscr E(\Psi(a)) +\mathscr E (\Psi(a))- \mathscr E(x_c) |^{1-\theta} \\
&\ge  |\mathscr E(x) - \mathscr E(x_c)|^{1-\theta} - C_3 \| \mathscr M(x)\|_{0,\alpha}^{2(1-\theta)}.
\end{split}
\end{equation}
Since $2(1-\theta) \ge 1$, (\ref{S L inequality}) is established for some $C_2 >0$ and for all $x\in \mathscr W_0$. 
\end{proof}

\subsection{Proof of Theorem \ref{finitely many entropy value on X_Lambda}}  
The following lemma is first proved in \cite{CM} (Lemma 7.10 therein) when $\Sigma$ is an $n$-dimensional self-shrinking embedded hypersurface in $\mathbb R^{n+1}$ with polynomial growth. Since a branched conformal immersion is immersed away from finitely many points, the exact same proof holds for compact branched conformally immersed self-shrinkers in $\mathbb R^m, m\geq 3$. For the reader's convenience, we sketch the proof of Lemma \ref{lambda = F_0,1} in the appendix. Note that the $\mathcal F$-functional (\ref{definition of mathcal F}) and the entropy (\ref{definition of entropy}) are also defined for branched immersions of compact surfaces.

\begin{lem} \label{lambda = F_0,1}
Let $F: \Sigma \to \mathbb R^m,m\geq 3$, be a compact branched conformally immersed self-shrinking surface. Then the entropy $\lambda$ defined in (\ref{definition of entropy}) is maximized at $(x_0, t_0) = (0,1)$. That is, 
\begin{equation}
\lambda(F) = \frac{1}{4\pi}\int_\Sigma e^{-\frac{|F|^2}{4}} d\mu. 
\end{equation} 
\end{lem}
Note that if $(F,\tau)$ is a critical point of $\mathscr E$, then $F$ is a branched conformally immersed self-shrinking surface. Conformality of $F$ then implies $|DF|^2_\tau d\mu_\tau = 2d\mu$, where 
$d\mu$ is the area element of the metric induced by $F$ away from the branch points. Together with (\ref{definition of mathscr E}) and Lemma \ref{lambda = F_0,1},
\begin{equation} \label{E = 4 pi lambda}
\mathscr E(F, \tau) = \int_\mathbb{T} e^{-\frac{|F|^2}{4}} d\mu = 4\pi \lambda(F). 
\end{equation}

Now we proceed to prove Theorem \ref{finitely many entropy value on X_Lambda}.

\begin{proof}
Assume the theorem is false. Then there is a sequence $\{F_n\} \in \mathfrak X_\Lambda$ with $\lambda(F_i)\neq \lambda (F_j)$ for all $i\neq j$. Let $g_n = F^*_n\langle\cdot, \cdot\rangle$ and let $g_{\tau_n}$ be the Riemannian metric on $\mathbb T$ which is of the form (\ref{definition of g_tau}) and is conformal to $g_n$. By Theorem \ref{compactness for branch self-shrinkers}, there is $F\in \mathfrak X_C$ and $\tau \in \mathbb H$ so that $F_n$ converges smoothly to $F$ and $\tau_n \to \tau$. Thus $$\| (F_n, \tau_n)- (F, \tau)\|_{2, \alpha} \to 0\,\,\,\,\mbox{as $n\to \infty$}.
$$
 From Proposition \ref{harmonic and conformal iff critical pt of E} and (\ref{E = 4 pi lambda}) and by setting $x_c = (F, \tau)$ in (\ref{S L inequality}),  we have $\lambda(F_i) = \lambda(F)$ for all $i$ large enough, since $\mathscr M(F_n, \tau_n) = 0$ for all $n$. That leads to a contradiction. Thus the theorem is proved.
\end{proof}


\section{Piecewise Lagrangian mean curvature flows}
In this section, we extend the definition of the piecewise MCF in \cite{CM} to Lagrangian MCF for torus in $\mathbb R^4$ and construct a piecewise Lagrangian MCF for a Lagrangian immersed torus $F: \mathbb T\to \mathbb R^4$. 

\begin{dfn}\label{definition of piecewise MCF}
Let $F: L\to\mathbb R^4$ be a Lagrangian immersion, where $L$ is a compact surface. A piecewise Lagrangian MCF with initial condition $F$ is a finite collection of smooth Lagrangian MCFs
$$F_t^i :  L  \to \mathbb R^4$$
defined on $[t_i, t_{i+1}]$, $i=0, 1, \cdots, k-1$, where $0=t_0 < t_1 < \cdots < t_{k-1} < t_k < \infty$  so that:
\begin{enumerate}
\item $F_0^0 =F$,
\item $\mu(F_{t_{i+1}}^{i+1}) = \mu(F_{t_{i+1}}^{i})$,
\item $\lambda( F_{t_{i+1}}^{i+1}) < \lambda (F_{t_{i+1}}^i)$,
\item there is $\delta >0$ such that 
\begin{equation} \label{perturbations are delta close}
\| F^{i}_{t_{i+1}} - F^{i+1}_{t_{i+1}}\|_{C^0} \le \delta \sqrt{\mu(F^i_{t_{i+1}})}
\end{equation} 
for $i = 0, 1, 2, \cdots, k-2$.
\end{enumerate}
\end{dfn}

\begin{rem}
Note that if $k=1$, the piecewise MCF is just the usual smooth MCF. The above definition is interesting only if we can characterize the behaviour of the flow when $t\to t_k$.
\end{rem}

Let $\{ F_t : L\to \mathbb R^4\}$ be a smooth MCF defined on $[t_0, T_0)$, where $T_0<\infty$ and $L$ is a closed surface.  Assume that a so-called type I singularity develops at $T_0$, which means  $\sup_{ F_t(L)} \|A_t\| \to \infty$ as $t\to T_0$ and  there is a positive constant $C$ so that

\begin{equation} \label{type I singularity definition}
\max_{F_t(L)} |A_t|^2 \le \frac{C}{\sqrt{T_0-t}}
\end{equation}
for all $t<T_0$. Let  $t_n\to T_0$ and $q_n\in F_{t_n}(L)$ where $\max_{F_n(L)} |A_{t_n} |$ is attained,  and suppose $q_n\to q\in\mathbb R^4$. 
Consider the type I rescaling, which is the family of immersions $\widetilde F(\cdot, s)$, where $-\log T_0 \le s < \infty$ and
\begin{equation} \label{type I blow up}
 \widetilde F(\cdot, s) = \frac{1}{\sqrt{(T_0 - t)}} (F_t(x) - q), \ \ s(t) = -\log (T_0 - t).
\end{equation}
For any sequence $s_j \to \infty$, 
a subsequence of $\{\widetilde F (\cdot, s_j)\}$ converges locally smoothly to a self-shrinking immersion $\overline F : \Sigma \to \mathbb R^4$ (\cite{H2}). In this case, we say that {\it the type I singularity can be modelled by $\overline F$}. It is not known whether $\overline F$ is unique: If we choose another sequence $\tilde s_k$, $\{\widetilde F(\cdot, \tilde s_k)\}$ might converge to a different self-shrinker.  


\vspace{.2cm}

Now we prove Theorem \ref{Generic LMCF}. 

\begin{proof} Let $F : \mathbb T \to \mathbb R^4$ be a Lagrangian immersion. By \cite{Sm1}, there is a unique smooth Lagrangian MCF $\{ F_t\}$ which is defined on a maximal time interval $[0,T_0)$, where $T_0 < \infty$ as $\mathbb T$ is compact. 

If the singularity at $T_0$ is not a type I singularity that can be modelled by a compact self-shrinker with area no larger than $\Lambda$, then we set $k=0$ and no perturbation is performed.

Otherwise, the singularity at $T_0$ is of type I and it can be modelled by a compact self-shrinker with area no larger than $\Lambda$. In this case, the inequality (\ref{type I singularity definition}) is satisfied at a point $q \in \mathbb R^4$ at time $T_0$ for some positive constant $C$ and for all $t\in [0,T_0)$, and there is a sequence $s_j \to \infty$ such that $\widetilde F(\cdot, s_j)$ as in (\ref{type I blow up}) converges locally smoothly to a compact self-shrinker $\overline F$ with area no bigger than $\Lambda$. To be precise about the convergence,  we recall that Lemma 3.3, Corollary 3.2 and Proposition 2.3 in \cite{H2} hold for any codimension, and they guarantee that all $\widetilde{F}(\cdot,s_j)$ touch a fixed bounded region, the areas inside a ball $B(R)$ are bounded by $C(R)$ and the second fundamental forms and their derivatives of any order are bounded. Therefore, 
all the conditions in Theorem 1.3 in \cite{Breunung} are satisfied for the sequence $\{\widetilde F(\cdot, s_j)\}$, and the theorem asserts: by passing to a subsequence if necessary, there is a surface $\Sigma$ and an immersion $\overline F : \Sigma \to \mathbb R^4$ and a sequence of diffeomorphisms 
$$
\varphi_j: U^j\to \widetilde F(\cdot, s_j)^{-1}(B_j) \subset\mathbb T,
$$ 
where $B_j$ is the ball of radius $j$ in $\mathbb R^4$ centered at the origin, $U_j \subset\Sigma$ are open sets with $U_j \subset\subset U_{j+1}$ and $\Sigma = \bigcup_j U_j$, such that 
$$\|\widetilde F(\cdot, s_j) \circ \varphi_j - \overline F\|_{C^0(U_j)}\to 0$$
and $\widetilde F(\cdot, s_j) \circ \varphi$ converges to $\overline F$ locally smoothly. 
In our situation, we have assumed that $\Sigma$ is compact (as we are dealing with singularity that can be modelled by compact shrinkers). Hence $\Sigma = U_k$ for all $k$ large and thus $\varphi_k $ are diffeomorphisms from $\Sigma$ to $\mathbb T$, since the torus is connected. To simplify notations, we write $\Sigma = \mathbb T$. The diffeomorphisms $\varphi_j : \mathbb T \to \mathbb T$ have the property that 
\begin{equation} \label{type I recscaling converges to C^k}
\| \widetilde F(\cdot, s_j)\circ \varphi_j - \overline F\|_{C^k(\mathbb T)} \to 0
\end{equation}
for all $k=0,1,2,\cdots$. Since each $\{F_t\}$ is Lagrangian, the sequence of blowups $\widetilde F(\cdot, s_j)$ are also Lagrangian for all $j$. The above convergence implies that $\overline F$ is Lagrangian, hence, $\overline F \in \mathfrak X_\Lambda$.

Since the entropy $\lambda$ (\ref{definition of entropy}) is translation and scaling invariant,
\begin{equation} \label{entropy translation and scaling invariant}
\lambda ( \widetilde F(\cdot, s(t))) = \lambda (F_t).
\end{equation}
Furthermore, by the definition of $\mathscr F_{x_0, t_0}$ in (\ref{definition of mathcal F}), we see 
\begin{equation} \label{entropy invariant under diff}
\lambda (\widetilde F (\cdot, s_j) \circ \varphi_j) = \lambda ( \widetilde F(\cdot, s_j)).
\end{equation}
Since $\mathscr F_{0,1}$ (see (\ref{definition of mathcal F})) is continuous with respect to the $C^1$-topology, there is a sequence $d_j$ of positive numbers so that $d_j \to 0$ as $j \to \infty$ and
\[
\mathscr F_{0,1} (\widetilde F(\cdot, s_j) \circ \varphi_j) \ge \mathscr F_{0,1} (\overline F) - d_j.
\]
By definition of $\lambda$ and Lemma \ref{lambda = F_0,1}, since $\overline F$ is a self-shrinker, from the above we have 
\begin{equation} \label{entropy bound for F( . , s_j)}
\lambda (\widetilde F (\cdot, s_j) \circ \varphi_j) \ge \lambda (\overline F ) - d_j. 
\end{equation}
As $\lambda$ is non-increasing along the MCF, $\lambda (\widetilde F(\cdot, s_j))$ is non-increasing in $j$ by (\ref{entropy translation and scaling invariant}). Together with (\ref{entropy invariant under diff}) and (\ref{entropy bound for F( . , s_j)}), we conclude
\[
\lambda (\widetilde F(\cdot, s_j ) ) \ge \lim_{j\to \infty} \lambda (\widetilde F(\cdot, s_j)) \ge \lambda (\overline F).
\]

Fix $\delta >0$. Let 
$$\delta_1 = \frac{\delta \sqrt{\mu(\overline F)}}{6},\ \ \ \delta_2 = \min\left\{ \frac 12, \frac{\delta_1 }{\|\overline F\|_{C^0}+\delta_1}\right\}.$$
Using (\ref{type I recscaling converges to C^k}), for all $k\ge 1$, there is $j_0$ so that 
\begin{equation} \label{F^R close to tilde F}
\| \widetilde F(\cdot, s_{j_0}) \circ \varphi_{j_0} - \overline F\|_{C^k} < \delta_1, 
\end{equation}
and 
\begin{equation} \label{mu (F^R) comparable to mu(tilde F)}
\left| \frac{\mu(\widetilde F(\cdot, s_{j_0}))}{\mu(\overline F)} -1\right|\le \delta_2 .
\end{equation}
By Theorem \ref{Theorem of Li and Zhang}, $\overline F$ is Lagrangian $\mathcal F$-unstable. Then by Theorem \ref{entropy unstable}, there is a Lagrangian immersion $\widehat F : \mathbb T \to \mathbb R^4$ which satisfies 
\begin{equation} \label{hat F close to tilde F}
\| \widehat F - \overline F\|_{C^2} <\delta_1,
\end{equation}
\begin{equation} \label{mu(tilde F) comparable to mu(hat F)}
  \left| \frac{\mu(\overline F)}{\mu(\widehat F)} -1\right|\le \delta_2
\end{equation}
and
\begin{equation}
\lambda (\widehat F) < \lambda (\overline F). 
\end{equation}

Now we define the first part of the piecewise Lagrangian MCF:
\begin{enumerate}
\item[(i)] The first piece of Lagrangian MCF is just $F^0_{t} := F_t$, where $t\in [0,t_1]$ and $t_1<T_0$ is such that $s(t_1) = s_{j_0}$. 
\item[(ii)] Define the first perturbation $F^1_{t_1}$ at time $t_1$ as 
\begin{equation} \label{definition of the first perturbation}
F^1_{t_1} = \sqrt{T_0 - t_1} ( \kappa \widehat F) \circ \varphi_{j_0}^{-1} + q.
\end{equation}
\end{enumerate}
where the dilation factor 
$$\kappa = \sqrt{\frac{\mu(\widetilde F (\cdot, s_{j_0})}{\mu(\widehat F)}}.$$
The constant $\kappa$ is chosen so that
\begin{equation} \label{mu(F'') = mu(F^R)}
\mu(\kappa \widehat F) = \mu(\widetilde F(\cdot, s_{j_0})).
\end{equation}

We check now that (2)-(4) in definition \ref{definition of piecewise MCF} are satisfied with $i=0$. First note that (2) follows from (\ref{mu(F'') = mu(F^R)}) and the definition of $F^0_{t_1}$ and $F^1_{t_1}$. To prove (3), since the entropy (\ref{definition of entropy}) is scaling and translation invariant, using $\lambda(\overline F)> \lambda(\widehat F)$ we obtain
\[
\lambda (F^0_{t_1}) = \lambda(\widetilde F(\cdot, s_{j_0}) \ge \lambda (\overline F)> \lambda (\widehat F) = \lambda (F^1_{t_1}).
\]
Thus (3) is also shown. Lastly, we show that (\ref{perturbations are delta close}) is satisfied with $i=0$. From (\ref{type I blow up}) and (\ref{definition of the first perturbation}), we have
\[ \|F^0_{t_1} - F^1_{t_1}\|_{C^0} = \sqrt{T_0 - t_1} \| \widetilde F(\cdot, s_{j_0})\circ \varphi_{j_0}- \kappa \widehat F\|_{C^0} .
\]
Note that (\ref{mu (F^R) comparable to mu(tilde F)}) and (\ref{mu(tilde F) comparable to mu(hat F)}) imply
\begin{equation} \label{mu(F^R) comparable to mu(hat F)}
\left| \kappa -1\right| \le B. 
\end{equation}
Together with (\ref{hat F close to tilde F}), (\ref{F^R close to tilde F}), the definition of $\delta_2$, we have 
\[
\begin{split}
\| \widetilde F(\cdot, s_{j_0})\circ \varphi_{j_0}- \kappa \widehat F\|_{C^0} &\le  \| \widetilde F(\cdot, s_{j_0})\circ \varphi_{j_0}- \overline F\|_{C^0} +\|\overline F - \widehat F\|_{C^0} +\| (1-\kappa)\widehat F\|_{C^0} \\
&\le 2\delta_1 + \delta_2(\delta_1+ \| \tilde F\|_{C^0}) \\
&\le 3\delta_1,
\end{split}
\]
where we used the simple estimate
$$
\|\widehat F\|_{C^0} \le \| \widehat F - \overline F \|_{C^0} + \|\overline F\|_{C^0}.
$$ 
Thus we have 
\[\begin{split}
\|F^0_{t_1} - F^1_{t_1}\|_{C^0} &< 3\delta_1\sqrt{T_0 - t_1} \\
&= 3\delta_1 \sqrt{\frac{\mu(F^0_{t_1})}{\mu(\widetilde F(\cdot, s_{j_0}))}} \\
&\le  \frac 12 \delta \sqrt{\mu(F^0_{t_1})} \sqrt{\frac{\mu(\overline F)}{\mu(\widetilde F(\cdot, s_{j_0}))}} \\
&\le \delta \sqrt{\mu(F^0_{t_1})},
\end{split}\]
where in the last step we used $\delta_2 \le \frac 12$. Thus (\ref{perturbations are delta close}) is shown and this finishes the construction of the first piece of the piecewise Lagrangian MCF. 

Using $F^1_{t_1}$ as initial condition, there is another family $\{ F_t : t\in [t_1, T_1)\}$ of smooth Lagrangian MCF with $F_{t_1} = F^1_{t_1}$. Again, if the condition in Theorem \ref{Generic LMCF} is satisfied at the singular time $T_1$ (that is, the singularity at $T_1$ is not of type I which can be modelled by a compact self-shrinker of area $\le \Lambda$), then we set $k=1$, $t_2 = T_1$ and we are done. If not, we carry out exactly the same procedure as above. Thus we have a Lagrangian self-shrinking torus $\overline F_1 \in \mathfrak X_\Lambda$, some time $t_2<T_1$ and another Lagrangian immersion $F^2_{t_2}$ so that 
$$\lambda (F^1_{t_2} )\ge \lambda (\overline F^1) > \lambda (F^2_{t_2}),$$
$$\mu(F^1_{t_2}) = \mu(F^2_{t_2})$$
 and 
 $$\|F^1_{t_2} -  F^2_{t_2}\|_{C^0} <\delta \sqrt{\mu(F^1_{t_2})}.
 $$
  Then, again, we apply the smooth Lagrangian MCF to $F^2_{t_2}$. Note that the above procedure must stop: Indeed, by Theorem \ref{finitely many entropy value on X_Lambda}, the image of $\lambda : \mathfrak X_\Lambda \to \mathbb R$ is finite. Moreover, from the above construction, each perturbation is chosen so that the entropy value is strictly less then one of the element in $\lambda(\mathfrak X_\Lambda)$. Since $\lambda$ is non-increasing along the usual MCF, the above procedure must terminate after $k$ steps for some $k\le |\lambda(\mathfrak X_\Lambda)|$. This implies that at $t_k$, the piecewise Lagrangian MCF do not encounter a type I singularity which can be modelled by a compact self-shrinker with area less than or equals to $\Lambda$.

To prove the last statement of Theorem \ref{Generic LMCF}, recall that the Maslov class of a Lagrangian immersion is given by $2[H] \in H^1(\mathbb T, \mathbb Z)$, where $H$ is the mean curvature form and $[H]$ is an integral class as
\begin{equation}\label{H = d theta}
H = d\theta,
\end{equation}
where $\theta : \mathbb T \to \mathbb S^1$ is the Lagrangian angle of the immersion $F : \mathbb T \to \mathbb R^4$ \cite{HarveyLawson}. When $\{F_t\}$ is a smooth Lagrangian MCF, $[H_t]$ is invariant as $[H_t]$ is an integral class and $H_t$ is smooth in $t$. This fact can also be checked using the evolution of $H$ under the Lagrangian MCF, see Theorem 2.9 in \cite{Sm1}. From (\ref{H = d theta}) it is also clear that the Maslov class is invariant under translation and scaling of the immersion. Thus when there is a type I singularity and $\overline F : \mathbb T \to \mathbb R^4$ is a compact Lagrangian self-shrinker which models the singularity, then $[H_{\overline F}] = [H_t]$. Lastly, we recall that in Theorem \ref{entropy unstable} the perturbation $\widehat F$ is defined using a closed 1-form on $\mathbb T$. Hence we also have $[H_{\overline F}] = [H_{F^1_{t_1}}]$. Thus the Maslov class is preserved when we perturb the Lagrangian immersion in constructing the piecewise Lagrangian MCF. This completes the proof of Theorem \ref{Generic LMCF}.
\end{proof}

\subsection{Generalization to Lagrangian immersion of higher genus surfaces}
Theorem \ref{Generic LMCF} can be extended to genus $g> 1$ if we impose further assumptions on the singularity. Let $c_1, c_2>0$ and consider the set $\mathfrak X^{\text{imm}}_{g, c_1, c_2}$ of all Lagrangian self-shrinking immersions $\overline F : \Sigma_g \to \mathbb R^4$ with area $\le c_1$ and the second fundamental form satisfying $\max_{\overline{F}(\Sigma_g)}|A|\le c_2$, where $\Sigma_g$ is a closed  orientable surface of genus $g$ with $g > 1$. Using \eqref{self-shrinking eqn}, there are constants  $C(k,c_1,c_2) >0$ that depend on $c_1,c_2,k$, such that
$$ 
\max_{\overline{F}(\Sigma_g)}|\nabla ^k A| \le C(k,c_1,c_2)
$$
for all $\overline F\in \mathfrak X^{\text{imm}}_{g, c_1, c_2}$. Thus we can apply Theorem 1.3 in \cite{Breunung} to conclude that $\mathfrak X^{\text{imm}}_{g, c_1, c_2}$ is compact in the $C^2$-topology, in particular, all sequential limits are unbranched. Unbranchedness of any limiting surface guarantees existence of nearby Lagrangian immersions by the Lagrangian neighbourhood theorem. By Theorem \ref{Theorem of Li and Zhang} and Theorem \ref{entropy unstable} again, the Lagrangian self-shrinkers in $\mathfrak X^{\text{imm}}_{g, c_1, c_2}$ are Lagrangian entropy unstable. It follows that all $\overline F\in \mathfrak X^{\text{imm}}_{g, c_1, c_2}$ are Lagrangian entropy unstable.  With these facts,  the proof of the following proposition is identical to that of Corollary 8.4 in \cite{CM} and is omitted here. 

\begin{prop} \label{lowering entropy on mathfrak X_g}
Let $\delta >0$. Then there is a positive constant $c$ depending only on $\delta$ such that for any Lagrangian self-shrinker $\overline F \in \mathfrak X^{\emph{imm}}_{g, c_1, c_2}$, there is a Lagrangian immersion $\widehat F : \Sigma_g \to \mathbb R^4$ so that $\| \widehat F - \overline F\|_{C^0} < \delta\sqrt{\mu(\overline F)}$ and $\lambda(\widehat F) < \lambda (\overline F) - c$.
\end{prop}

\begin{rem}
For genus $>1$, without assuming uniform boundedness of the second fundamental forms, the compactness result in Theorem \ref{compactness for branch self-shrinkers} is not enough to conclude Proposition \ref{lowering entropy on mathfrak X_g} due to the assumption on the conformal structures there. 
\end{rem}

Using Proposition \ref{lowering entropy on mathfrak X_g}, we can define a piecewise Lagrangian MCF for a Lagrangian immersion $F: \Sigma_g \to \mathbb R^4$, as we did in the genus 1 case. After each perturbation, the entropy decreases by a fixed amount $c>0$ (Note that this $c$ might depend on $\delta$). Since the entropy is always is positive number, we conclude that the process must terminate in finite time and we have the following

\begin{thm}\label{g>1}
Let $F : \Sigma_g \to \mathbb R^4$ be a Lagrangian immersion and $\delta >0$ be given. Then there exists a piecewise Lagrangian MCF $\{F^i_t : i = 0, 1, \cdots, k-1\}$ with initial condition $F$, such that the singularity at time $t_{k}$ is not a type I singularity which can be modelled by a self-shrinker in $\mathfrak X^{\emph{imm}}_{g, c_1, c_2}$.  Moreover, we have the estimates $\| F^i_{t_i} - F^{i+1}_{t_i}\|_{C^0} <\delta\sqrt{\mu(F^i_{t_i})}$ and the Maslov class of each immersion is invariant along the flow.
\end{thm}

\section{Appendix}
\subsection{Proof of Lemma \ref{lambda = F_0,1}} Let $\Sigma$ be a compact surface without boundary and let $F:\Sigma \to \mathbb R^4$ be a branched conformal self-shrinker. Define the operator $\mathcal L_s$ by
\begin{equation}
\mathcal L_s u = \Delta u - \frac{1}{2t_s} \il (x-x_s)^\top , \nabla u\li  = e^{\frac{|x-x_s|^2}{4t_s} }\text{div} (e^{-\frac{|x-x_s|^2}{4t_s}} \nabla u).
\end{equation}
Here $(x_s, t_s) \in \mathbb R^4 \times \mathbb R_{>0}$, $\nabla$, $\text{div}$ and $\Delta$ are taken with respect to the pullback metric $F^*\langle\cdot, \cdot\rangle$ and $u, v$ are functions on $\mathbb R^4$. Note that $\mathcal L_s$ is defined away from the set of branch points $B$. As in \cite{CM}, we use the square bracket $[\cdot]_s$ to denote
\begin{equation}
[f]_s = \frac{1}{4\pi t_s}\int_\Sigma f e^{-\frac{|x-x_s|^2}{4t_s}} d\mu
\end{equation}

\begin{lem} We have 
\begin{equation} \label{L self adjoint}
[ u \mathcal L_s v]_s   =  -[\langle \nabla u, \nabla v\rangle ]_s . 
\end{equation}
\end{lem}

\begin{proof}
Let $B = \{x_1, \cdots, x_n\}$. Let $\epsilon>0$ be small and $B_i(\epsilon)$ be an $\epsilon$-ball in $\Sigma$ with center $x_i$, so that $B_i(\epsilon)\cap B_j(\epsilon) = \emptyset $ if $i\neq j$. Then

\begin{equation}\begin{split}
[u \mathcal L v ]_s &= \frac{1}{4\pi t_s}\int_\Sigma u \,\text{div} \left( e^{-\frac{|x-x_s|^2}{4t_s}} \nabla v\right) d\mu \\
&= \lim_{\epsilon \to 0} \frac{1}{4\pi t_s} \int_{\Sigma \setminus \bigcup B_i(\epsilon)} u \,\text{div} \left( e^{-\frac{|x-x_s|^2}{4t_s}} \nabla v\right) d\mu \\
&= \lim_{\epsilon \to 0}\frac{1}{4\pi t_s} \left( \sum_i \int_{\partial B_i(\epsilon)} u \il \nabla v, n_i \li  e^{-\frac{|x-x_s|^2}{4t_s}} dl- \int_{\Sigma \setminus \bigcup B_i(\epsilon)} \il \nabla u, \nabla v\li e^{-\frac{|x-x_s|^2}{4t_s}} d\mu \right)\\
&= - [\langle \nabla u, \nabla v\rangle ]_s
\end{split}
\end{equation}
where $n_i$ is the unit outward normal along $\partial B_i(\epsilon)$.
\end{proof}

In particular, we have 
\begin{equation} \label{mathcal L_s self adjoint appendix}
[u\mathcal L_s v]_s = -[\langle \nabla u, \nabla v\rangle]_s = [v\mathcal L_s u]_s.
\end{equation}
Using (\ref{mathcal L_s self adjoint appendix}), exactly the same argument in \cite{CM}, pp. 786-788, shows that for all $y\in\mathbb R^4$ and $a\in \mathbb R$ if we set $(x_s, t_s) = (sy, 1+ as^2)$ and $g(s) =\mathcal F_{x_s, t_s} (F)$ then $g'(s) \le 0$ for all $s>0$ with $1+as^2 >0$. Thus $\mathcal F_{y, t} (F) \le \mathcal F_{0,1}(F)$ for all $(y, t) \in \mathbb R^4 \times \mathbb R_{>0}$ and thus Lemma \ref{lambda = F_0,1} is proved.

\subsection{Analyticity of $\mathscr E$ and $\mathscr M$} Next we show that both $\mathscr E$ and $\mathscr M$ defined in (\ref{definition of mathscr E}) and (\ref{definition of mathscr M}) are analytic. For the definition of continuous symmetric $n$-linear form and analytic function between Banach spaces, please refer to Chapter 4 in \cite{Deimling}. First we have 

\begin{lem} \label{product is analytic}
Let $X$, $Y$ $Z$ be Banach spaces, $\mathscr U$, $\mathscr V$ are open in $X$, $Y$ respectively, and $f : \mathscr U \to \mathbb R$, $g : \mathscr V \to Z$ are analytic at $x_0\in \mathscr U$, $y_0 \in\mathscr V$ respectively. Then the function 
$$h : \mathscr U \times \mathscr V \to Z,\ \ \ h (x, y) = f(x) g(y)$$ 
is analytic at $(x_0, y_0)$.
\end{lem}

\begin{proof}
Since $f, g$ are analytic at $x_0, y_0$ respectively, then 
\begin{equation} \label{power series expansion of f and g}
f(x_0 + h) = f(x_0) + \sum_{i=1}^\infty A_i(h^i), \ \ \ g(y_0 + k) = g(y_0) + \sum_{j=1}^\infty B_j (k^j)
\end{equation}
for all $\|h\|_X<\epsilon_1$ and $\|k\|_Y<\epsilon_2$, and $A_i, B_j$ are continuous multi-linear forms such that 
\begin{equation} \label{bounds for A_i and B_j}
\sum_{i=1}^\infty \|A_i\| \epsilon_1^i <+\infty \,\,\,\mbox{and} \,\,\,\, \sum_{j=1}^\infty \|B_j\| \epsilon_2^j <+\infty.
\end{equation}
The absolute convergence of (\ref{power series expansion of f and g}) implies that 
\begin{equation} \label{expansion of product}
h(x_0+h, y_0+k) = f(x_0) g(y_0) +\sum_{n=1}^\infty C_n(h,k), 
\end{equation} 
for all $(h, k)$ such that $\|h\|_X <\epsilon_1, \|k\|_Y<\epsilon_2$, where 
\begin{equation}
C_n (h,k) = \sum_{i = 0}^n A_i(h^i) B_{n-i}(k^{n-i}).
\end{equation}
Let $\epsilon = \frac 12 \min\{ \epsilon_1, \epsilon_2\}$. Then by definition of $\|C_n\|$ and $\epsilon$, one has
\[
\begin{split}
\|C_n\| \epsilon^n &= \sup_{\|h\|+\|k\| = \epsilon} \|C_n(h,k)\|_Z \\
&\le  \sum_{i=0}^n\left( \| A_i\| \epsilon_1^i\right)\left( \| B_{n-i}\| \epsilon_2^{n-i}\right). 
\end{split}
\]
Thus 
\[
\sum_{n=1}^\infty \|C_n\| \epsilon^n \le\left(\sum_{i=1}^\infty \|A_i\| \epsilon_1^i\right)\left(\sum_{j=1}^\infty \|B_j\| \epsilon_2^j\right) <+\infty
\]
by (\ref{bounds for A_i and B_j}). Hence $h$ is also analytic at $(x_0, y_0)$. 
\end{proof}

\begin{prop} \label{mathscr E is analytic}
The mapping $\mathscr E:\mathscr U  \to \mathbb R$ in (\ref{definition of mathscr E}) is analytic.
\end{prop}

\begin{proof}
Using (\ref{definition of g_tau}), we have 
\begin{equation}
2\mathscr E(u, \tau) = \left(\tau^2_1/\tau_2+\tau_2\right) L_{11}(u) - (2\tau_1/\tau_2) L_{12}(u) +L_{22}(u),
\end{equation}
where 
\begin{equation}
L_{ij}(u) = \int_{\mathbb T} D_i u \cdot D_j u e^{-\frac{|u|^2}{4}} dxdy. 
\end{equation}
Since $\tau\mapsto (\tau_1^2/\tau _2) +\tau_2$ and $\tau \mapsto \tau_1/\tau_2$ are analytic, by Lemma \ref{product is analytic}, it suffices to check $L_{ij} : C^{2,\alpha} \to \mathbb R$ is analytic. But this is obvious, using the power series expansion of $e^{-\frac{|u|^2}{4}}$. 
\end{proof}

\begin{prop}
The mapping $\mathscr M :\mathscr U \mapsto \mathscr C^{0,\alpha}$ in (\ref{definition of mathscr M}) is analytic. 
\end{prop}

\begin{proof}
It suffices to show that both components in (\ref{expression of mathscr M}) are analytic. The second component $(u, \tau)\mapsto \nabla \mathscr E^u_\tau$ is analytic since $\mathscr E$ is analytic by Proposition \ref{mathscr E is analytic}, here we recall that $\nabla\mathscr E^u_\tau$ is the gradient of $\mathscr E(u,\tau)$ at $\tau$. Note that the first component can be written as 
\begin{equation} \label{second component of mathscr M in appendix}
(u, \tau)\mapsto -g^{ij}_\tau \left( D_{ij} u - (u\cdot D_j u) D_i u +\frac 14 (D_i u\cdot D_ju)u\right)
\end{equation}
Since $\tau \mapsto g^{ij}_\tau$ is analytic, the mapping in (\ref{second component of mathscr M in appendix}) is also analytic by Lemma \ref{product is analytic}. 
\end{proof}

\end{document}